\documentclass[12pt]{article}
\usepackage{amsmath}
\usepackage{amsthm}
\usepackage{amssymb}
\newtheorem{defin}{Definition}[section]
\newtheorem{prop}[defin]{Proposition}

\newtheorem{lemma}[defin]{Lemma}
\newtheorem{nota}[defin]{Notation}
\newtheoremstyle{recht}
  {\topsep}
  {\topsep}
  {\rmfamily}
  {}
  {\bfseries}
  {.}
  {.5em}
  {}
{\theoremstyle{recht}
\newtheorem{remark}[defin]{Remark}}
\newcommand{\C}{\mathbb{C}}
\newcommand{\vphi}{\varphi}
\newcommand{\N}{\mathbb{N}}
\newcommand{\ot}{\otimes}
\newcommand{\oh}{\otimes_h}
\newcommand{\od}{\odot}
\newcommand{\cL}{\mathcal{L}}
\newcommand{\Hfi}{H_\varphi}
\newcommand{\Hpsi}{H_\psi}
\newcommand{\lapsi}{\Lambda_\psi}
\newcommand{\pipsi}{\pi_\psi}
\newcommand{\ronda}{\mathcal{A}}
\newcommand{\Nfi}{\mathcal{N}_\varphi}
\newcommand{\Mfi}{\mathcal{M}_\varphi}
\newcommand{\Npsi}{\mathcal{N}_\psi}
\newcommand{\pifi}{\pi_\varphi}
\newcommand{\lafi}{\Lambda_\varphi}
\newcommand{\io}{\iota}
\newcommand{\eps}{\varepsilon}
\newcommand{\om}{\omega}
\newcommand{\omh}{\otimes_{mh}}
\begin{document}
\title{Hopf C$^*$-algebras}
\author{Stefaan Vaes\footnote{Research Assistant of the Fund for Scientific
Research - Flanders (Belgium)(F.W.O.)} \footnote{Supported by the F.W.O. Research Network WO.011.96N}
$\quad$ and $\quad$ Alfons Van Daele \\
Department of Mathematics \\ Katholieke Universiteit Leuven
\\ Celestijnenlaan 200B \\ B-3001 Heverlee \\ Belgium\\ email :
Stefaan.Vaes@wis.kuleuven.ac.be \\
email : Alfons.VanDaele@wis.kuleuven.ac.be
\\ fax : (32)-16 32 79 98}
\date{\textbf{July 1999}}
\maketitle
\begin{abstract}
\setlength{\parskip}{1ex plus 0.5ex minus 0.2ex}

\noindent In this paper we define and study Hopf C$^*$-algebras.  Roughly speaking, a
Hopf C$^*$-algebra is a C$^*$-algebra $A$ with a comultiplication $\phi :
A \rightarrow M(A \otimes A)$ such that the maps $a \otimes b
\mapsto \phi(a)(1 \otimes b)$ and $a \otimes b \mapsto  (a \otimes 1)\phi(b)$ have
their range in $A \otimes A$ and are injective after being extended to a
larger natural domain, the Haagerup tensor product $A \oh A$.
In a purely algebraic setting, these conditions on $\phi$ are
closely related with the existence of a counit and antipode.  In this
topological context, things turn out to be much more subtle, but nevertheless
one can show the existence of a suitable counit and antipode
under these conditions.

\noindent The basic example is the C$^*$-algebra $C_0(G)$ of continuous complex
functions tending to zero at infinity on a locally compact group where the
comultiplication is obtained by dualizing the group multiplication.  But also the
reduced group C$^*$-algebra $C^*_r(G)$ of a locally compact group with
the well known comultiplication falls in this category.  In fact all locally compact
quantum groups in the sense of \cite{KV1} and \cite{KV2} (such as the compact
and discrete ones) as well as most of the known examples are included.

\noindent This theory differs from other similar approaches in
that there is no Haar measure assumed.
\end{abstract}

\section*{Introduction}
Let $G$ be a locally compact space.  Consider the C$^*$-algebra $C_0(G)$ of
all complex continuous functions on $G$ tending to $0$ at infinity.  Any
abelian C$^*$-algebra has this form.  Moreover, the topological structure of
$G$ is completely determined by the C$^*$-algebra structure of $C_0(G)$.
This is why an arbitrary C$^*$-algebra is sometimes thought of as `the
continuous functions tending to $0$ at infinity on a pseudo space' (a quantum
space or a non-commutative space).  In this point of view, the study of
C$^*$-algebras is non-commutative topology.

Now, one can  consider extra structure on such a locally compact space,
translate this structure to the C$^*$-algebra and then formulate it on the
abstract level, for any C$^*$-algebra.  This is the quantization
procedure.  This paper deals with quantizing the group structure of a locally
compact group.

Before we can explain this, we need to recall some results in C$^*$-algebra
theory.  When $A$ is a C$^*$-algebra, we denote by $M(A)$ the multiplier
algebra of $A$.   A $^\ast$-homomorphism $\gamma : A
\rightarrow M(B)$ for C$^*$-algebras $A$ and $B$ is called {\it
non-degenerate} if an approximate identity in $A$ is mapped to a net strictly
converging to 1 in $M(B)$.
This is equivalent with the requirement that $\gamma(A)B$ is dense in $B$.
In Woronowicz' terminology one would call $\gamma$ a morphism from
$A$ to $B$.
An important property of such a non-degenerate $^\ast$-homomorphism is that it has
a unique extension to a unital $^\ast$-homomorphism from $M(A)$ to $M(B)$.
This extension is still denoted by $\gamma$ and it is defined by
$$\gamma(x) \gamma(a)b = \gamma(xa)b$$
whenever $x \in M(A), a \in A$ and $b \in B$.  The extension is strictly
continuous on bounded subsets of $M(A)$.

Now, let $G$ be a locally compact group.  First, we will explain how the group
structure is translated to the C$^*$-algebra $C_0(G)$.  Take $f
\in
C_0(G)$ and dualize the product in $G$ to get a function $\phi(f)$ of two
variables defined by $\phi(f)(p,q) = f(pq)$ for any $p,q \in G$.  This function
is bounded and continuous.  We clearly get a $^\ast$-homomorphism
$\phi$ from $C_0(G)$ to $C_b (G \times G)$, the C$^*$-algebra of bounded,
continuous complex functions on $G \times G$.  Now, there is a natural
identification of $C_b(G \times G)$ with the multiplier algebra $M(C_0(G
\times G))$ and further of $C_0(G \times G)$ with $C_0(G) \otimes C_0(G)$, the
(minimal) C$^*$-tensor product of $C_0(G)$ with itself.  Therefore, we may
consider $\phi : C_0(G) \rightarrow M(C_0 (G) \otimes C_0(G))$.  This
$^\ast$-homomorphism is non-degenerate.  The associativity of the product in
$G$ gives rise to the coassociativity law $(\phi \otimes \iota)\phi = (\iota
\otimes \phi) \phi$ of $\phi$.  Here $\iota$ denotes the identity map.  And
$\phi \otimes \iota$ and $\iota \otimes \phi$ are the unique extensions to
$M(C_0(G) \otimes C_0(G))$ of the obvious maps on $C_0(G) \otimes C_0(G)$.

In definition~\ref{21}, we will generalize this notion and define a comultiplication
on a C$^*$-algebra $A$ as a non-degenerate $^\ast$-homomorphism $\phi : A
\rightarrow M(A \otimes A)$ satisfying coassociativity: $(\phi \otimes
\iota)\phi = (\iota \otimes \phi)\phi$.

It is important to observe that, conversely, any comultiplication on an abelian
C$^*$-algebra $A$ comes from a continuous associative multiplication of the underlying
locally compact space.  Therefore we will think of a pair $(A,\phi)$ of a
C$^*$-algebra $A$ with a comultiplication $\phi$ as a {\it locally compact
quantum semi-group}.

Next, consider the identity $e$ of the group.  This gives rise to a
$^\ast$-homomorphism $\varepsilon : C_0(G) \rightarrow \C$ given by $\varepsilon
(f) = f(e)$.  The property $pe = ep = p$ for all $p \in G$ is expressed in
terms of $\varepsilon$ by
$$(\varepsilon \otimes \iota)\phi(f) = (\iota \otimes \varepsilon)\phi(f) = f$$
for all $f \in C_0(G)$.  This naturally leads us to the notion of a counit for a
general pair $(A, \phi)$.  A possible definition of a
counit is a $^\ast$-homomorphism $\varepsilon : A \rightarrow  \C$ satisfying
$$(\varepsilon \otimes \iota) \phi (x) = (\iota \otimes \varepsilon) \phi (x) = x$$
for all $x \in A$.  Such a notion makes sense but unfortunately turns out to be
too restrictive.  The typical example to illustrate this comes from a locally
compact non-amenable group $G$.  Take for $A$ the reduced group
C$^*$-algebra $C^\ast_r(G)$.  This carries a natural comultiplication $\phi$
characterized by $\phi(\lambda_p) = \lambda_p \otimes \lambda_p$ where
$\lambda_p$ is left translation by $p$ on $L^2(G)$, considered as a
multiplier of $C^\ast_r(G)$.  If $\varepsilon$ is a counit as above, it should
verify $\varepsilon(\lambda_p) = 1$ for all $p$.  However, as $G$ is not amenable,
the trivial representation of $G$ is not weakly contained in the regular
representation and so there is no such bounded $^\ast$-homomorphism on
$C^\ast_r(G)$ (see e.g.\ \cite{JPP}, Theorem~8.9).  This problem could be overcome by taking the full group
C$^*$-algebra $C^\ast(G)$.  This will create a new problem~: the natural
candidate for the Haar measure will no longer be faithful.  This turns out to
be more serious.

Therefore, we will work with another notion of a counit (see proposition~\ref{38}),
possibly not everywhere defined (and unbounded).  A $^\ast$-homomorphism
$\varepsilon : A \rightarrow  \cal \C$ satisfying $(\varepsilon \otimes \iota) \phi
(x)
= (\iota \otimes \varepsilon) \phi (x) = x$ for all $x \in A$ will be called a {\it
bounded} counit.

Finally, consider the inverse in the group. We will see that this is even more
problematic.  Define the map $S : C_0(G) \rightarrow C_0(G)$ by $(Sf)(p) =
f(p^{-1})$ for all $p \in G$ : this map is called the antipode.  Here it is a
$^\ast$-homomorphism.  The properties $pp^{-1} = e$ and $p^{-1}p = e$ for all
$p \in G$ are translated in terms of the comultiplication $\phi$, the counit
$\varepsilon$ and the antipode $S$ on the C$^*$-algebra level into the formulas
\begin{align*}
m(\iota \otimes S) \phi (f) & = \varepsilon(f) 1 \\
m(S \otimes \iota) \phi (f) & = \varepsilon(f) 1
\end{align*}
for all $f \in C_0(G)$.  Here $m$ denotes the multiplication map, defined from
$C_b (G \times G)$ to $C_b(G)$ by $m(h)(p) = h(p,p)$ when $h \in C_b(G \times
G)$.  Observe that $m(h) = fg$ when $h = f \otimes g$ as then
$$m(h)(p) = h(p,p) = f(p)g(p) = (fg)(p)$$
for all $p \in G$.

If we want to formulate all of this for a general C$^*$-algebra $A$ with a
comultiplication $\phi$, we run into all sorts of problems.  First, we need a
counit $\varepsilon$ and this was already a problem. Secondly, as we
know from general Hopf algebra theory, as soon as we work with non-abelian
algebras, the antipode $S$ should be an anti-homomorphism, rather than a
homomorphism.  Furthermore, in the interesting examples, the candidate for
the
antipode will not be a bounded everywhere defined map.  This looks similar to
the counit being unbounded, but it is a different problem.  In the example
$C^\ast_r(G)$ with $G$ non-amenable $\varepsilon$ is unbounded while $S$ is
bounded.  The unboundedness of $S$ is related with the property that $S$ also
need not be a $^\ast$-map.  Now, we are using $S \otimes \iota$ and $\iota
\otimes S$ for characterizing the antipode and it is not clear how to give a meaning
to these maps when $S$ is unbounded.  But even when $S$ were a
$^\ast$-anti-homomorphism, it is not obvious how to define this on the
C$^*$-tensor product $A \otimes A$ and so also not on $M(A \otimes A)$.
Finally the last problem that we are facing here is the multiplication map $m$.
It is in general defined on the algebraic tensor product $A \odot A$ in the
obvious way.  But generally it will not extend to $A \otimes A$ and so not to
$M(A \otimes A)$.  That this can be done is typical for the abelian case.

So, it should be clear from the above discussion that there are serious
problems with the concept of the antipode and (although less serious) with the
counit for a general pair of a C$^*$-algebra $A$ with a comultiplication
$\phi$.

This is probably the reason why, up to now, there was not a good notion of a
Hopf C$^*$-algebra.  In the literature, there existed essentially
two directions.  Some authors already used the name Hopf C$^*$-algebra for a
pair of a C$^*$-algebra with a comultiplication, see e.g.~\cite{B-S}, \cite{Val} and \cite{Ng}.  This is not
very appropriate because the antipode is necessary to distinguish a Hopf
algebra from a bialgebra.  In the other approach, the existence of an antipode
is assumed.  In the older literature (e.g.\ when defining C$^*$-Kac algebras,
see~\cite{E-V} and \cite{Val}), the map $S$ is taken to be a $^\ast$-anti-automorphism satisfying
$\phi(S(a)) = \sigma (S \otimes S) \phi (a)$ where $\sigma$ is the flip.  The
latter is a property of the antipode in a Hopf algebra but much weaker than the
usual requirement.  This shortcoming was dissolved by imposing an extra
condition involving the Haar measure.  Moreover, from recent examples, we know
that it is also too strong to require $S$ to be everywhere defined and to
satisfy $S(a^\ast) = S(a)^\ast$.  This last difficulty can be overcome by using
some kind of polar decomposition of the antipode.  But still, the basic axiom
for $S$ is formulated in connection with the Haar measure (see e.g.~\cite{Mas-Nak}).

In this paper we propose a solution to this problem.  We will start with a pair
of a C$^*$-algebra $A$ and a comultiplication $\phi$ on $A$.  And we will
impose extra, but natural conditions on $\phi$ so that it is possible to
construct a counit $\varepsilon$ and an antipode $S$ and to give a meaning to the
usual formulas

\begin{align*}
(\varepsilon \otimes \iota)\phi(a) & = (\iota \otimes \varepsilon)\phi(a) = a \\
m(S \otimes \iota)\phi(a) & = m (\iota \otimes S)\phi (a) = \varepsilon(a)1
\end{align*}
for all $a$ in a suitable subspace of $A$.

The approach will heavily rely on the Haagerup tensor product $A \otimes_h A$.
We will recall in section 1 the definition and the basic properties.  The main
reason for using the Haagerup tensor product is that it is a natural domain for
the multiplication map.  Indeed, the multiplication $m : A \odot A
\rightarrow A$ defined on the algebraic tensor product, is a contraction for
the Haagerup norm and hence extends by continuity to a contraction from the
completion $A \otimes_h A$.  In fact, we will need to work with some larger
space, denoted by $A \otimes_{mh} A$.  Elements $x$ in this space have the
property that $(a \otimes 1)x$ and $x(1 \otimes a)$ belong to $A  \otimes_h A$
for all $a \in A$.  This will be true essentially by definition.  This larger
space will also be introduced in section 1, where we will give its properties
needed further in the paper.  The basic point is that multiplication is now a
contraction from this large space $A \otimes_{mh} A$ to $M(A)$.

The conditions on $\phi$ are formulated in terms of the maps $T_1$ and $T_2$ defined
by
\begin{align*}
T_1 (a \otimes b) & = \phi(a)(1 \otimes b) \\
T_2 (a \otimes b) & = (a \otimes 1) \phi(b).
\end{align*}
These maps are first defined on the algebraic tensor product $A \odot A$ and
have range in $A \otimes A$ (by one of the assumptions on $\phi$).  Now, it
should be no surprise that these maps are continuous for the Haagerup norm on
$A \odot A$ and hence extend to contractions from $A \otimes_h A$ to $A
\otimes A$.  The extra requirement is that these extensions are injective.

Recall that in the case of a Hopf algebra (or more generally, a multiplier Hopf
algebra), these maps $T_1$ and $T_2$ are bijective from the algebraic tensor
product of the underlying algebra to itself.  In fact,
this is essentially the property characterizing Hopf algebras among bialgebras.
We obviously have a topological analogue of this condition here.
And we will see further in the paper why it is precisely the Haagerup tensor
product that is natural to use.

The conditions that we will impose will obviously be fulfilled in our
motivating example coming from a locally compact group.  It will also be
possible to verify them in many, if not all of the existing examples (see
section 5).

Much of this paper is inspired by the work on multiplier Hopf algebras.  In
fact, what we do here is, to a great extend, consider the different steps in
the original development of multiplier Hopf algebras (see~\cite{VD2}) and try to
translate this to the C$^*$-framework.  Therefore, many formulas in this
paper are inspired by their analogues in multiplier Hopf algebras.  We will
sometimes refer to these and use them as a motivation or to illustrate the new
formulas here.

However, it is not required to know the work on multiplier Hopf algebras.  It
is even not necessary to know Hopf algebras (although this will of course help
in understanding the paper).  What we do require is basic C$^*$-algebra
theory.

{\it The paper is organized as follows.}  In section 1 we start by recalling the
notion of the Haagerup tensor product and its basic properties.  We also define
a larger space and prove the necessary properties needed for the rest of the
paper.  In section 2 we recall the notion of a comultiplication $\phi$ on a
C$^*$-algebra $A$ and we explain the extra condition we will need to develop
the theory of Hopf C$^*$-algebras.  In section 3 we prove the existence of
the counit and the antipode under the given condition.  It should be said
however that we can obtain this from two directions, left and right.  In
section 4 we treat the construction of the counit and the antipode in a
\lq two-sided\rq\ way.  We will essentially get a domain on which the left and right
counits and antipodes will coincide.  Finally, in section 5 we will treat many
examples and thus justify our axioms and our definition of Hopf
C$^*$-algebras. Because we will be working with the Haagerup
tensor product, the theory of operator spaces will be in the
background throughout the paper although we never require more
then elementary results from that theory.

The paper is a slight modification of the first author's masters thesis~\cite{V}.
The main difference lies in the use of another type of generalized Haagerup
tensor product.  We will say more about this in section 2.  In principle, the
approach here is more general (although we still are not sure about this).  It
is somewhat more natural and in some sense simpler.

Let us now finish this introduction and fix some basic notations.  We will use
$(H,\Delta,\varepsilon,S)$ for a Hopf algebra (or a multiplier Hopf algebra) and
in that context, all the tensor products will be algebraic tensor products.  We
will however use $\odot$ for the algebraic tensor product when we are
dealing with topological spaces (C$^*$-algebras most of the time) and
$\otimes$ for completed tensor products.  In particular $A \otimes A$ will be
the minimal C$^*$-tensor product of the C$^*$-algebra $A$ with itself.  When
we write $\otimes_h$, we will also mean the completed Haagerup tensor product
(see section 1).  We will use $\iota$ to denote the identity maps.  We will
identify $A \otimes \C$ and $\C \otimes A$ with $A$ and so slice maps
$\iota \otimes \omega$ and $\omega \otimes \iota$ will go e.g.\ from $A \otimes
A$ to $A$ (for any $\omega \in A^\ast$, the norm dual of $A$).  As we
mentioned already before, $M(A)$ will denote the multiplier algebra of the
C$^*$-algebra $A$ and it will be considered both with its norm topology and
the strict topology making the maps $x \mapsto ax$ and $x \mapsto xa$
from $M(A)$ to $A$ continuous for all $a \in A$.  We have already mentioned the
extension property for non-degenerate $^\ast$-homomorphisms.  We will also need
to extend functionals in $A^\ast$.  It is important for this to notice that
elements $\omega \in A^\ast$ are always of the form $\omega = \omega_1 ( \, \cdot \, c) =
\omega_2 (d \, \cdot \, )$ for some $c,d \in A$ and $\omega_1, \omega_2
\in A^\ast$.  It follows that elements in $A^\ast$ have unique
strictly continuous extensions to $M(A)$. We will sometimes denote
$\om_1(\, \cdot \, c)$ by $c \om_1$ and $\om_2(d \, \cdot \,)$ by
$\om_2 d$. We will use the symbol $\phi$ to denote a
comultiplication on a C$^*$-algebra $A$ (see def.~\ref{21}). This
will be a non-degenerate $^*$-homomorphism from $A$ to $M(A \ot
A)$ (hence a morphism from $A$ to $A \ot A$ in Woronowicz' terminology)
satisfying $(\phi \ot \io) \phi = (\io \ot \phi) \phi$. Then
we use the notation $\phi^{(2)}=(\phi \ot \io) \phi = (\io \ot \phi)
\phi$. We will use sometimes the leg-numbering notation, as
well for $^*$-homomorphisms as for operators. So $\phi_{12}$ will
simply denote $\phi \ot \io$ while $\phi_{13}$ equals $\phi_{12}$
composed with the flip map between the second and the third leg.
Analogously one introduces $u_{12}$ and $u_{13}$ when $u$ is a
bounded operator or an element in some tensor product.

The basic references for C$^*$-algebras are~\cite{Ped} and \cite{Tak}.  For Hopf algebras we
refer to~\cite{Abe} and for multiplier Hopf algebras to~\cite{VD2}.

{\bf Acknowledgment.}  The second author likes to thank his colleagues at the
University of Fukuoka (where this work was started).  We are also grateful for
the hospitality we both experienced at the University of Trondheim (where this
work was completed).  Moreover we enjoy the pleasant working atmosphere in Leuven
(where most work on this paper was done).

\section{The multiplier Haagerup tensor product}
In this preliminary section we will first recall the definition of the
Haagerup tensor product and its properties we will need further in the paper.
For details, we refer to e.g.~\cite{B-P}, \cite{B-P2} and \cite{E-R}.

Let $A$ be a C$^*$-algebra.  Denote by $A \odot A$ the algebraic tensor product of
$A$ with itself.

\begin{defin} \label{11}
The Haagerup tensor product on $A \odot A$ is defined
by
$$\| x \|_h = \inf \Bigl\{ \Bigl\| \sum_{i=1}^n p_i p^*_i \Bigr\|^{1/2} \,  \Bigl\| \sum_{i=1}^n q^\ast_i q_i
\Bigr\|^{1/2} \, \Big| \,  n \in \N, x = \sum^n_{i=1} p_i \otimes q_i \Bigr\}.$$
\end{defin}
So, the infimum is taken over all the possible representations of $x$ in $A
\odot A$.  One can prove that this is indeed a norm.  The
completion of $A \odot A$ with respect to this norm will be denoted here by
$A \otimes_h A$.

Let now $I$ be an index set and suppose that $(p_i)_{i \in I}$ and
$(q_i)_{i \in I}$ are families of elements of $A$ such that
the nets $\sum_{i \in I} p_i p_i^*$ and $\sum_{i \in I}
q_i^*q_i$
are norm convergent in $A$. Then the net
$\sum_{i \in I} p_i \ot q_i$
will be norm convergent in $A \ot_h A$. This is quite easy to see,
because for $I_0 \subseteq I$ and $I_0$ finite we have
$$\Bigl\| \sum_{i \in I_0} p_i \ot q_i \Bigr\|_h^2 \leq
\Bigl\| \sum_{i \in I_0} p_i p_i^* \Bigr\| \, \Bigl\| \sum_{i \in I_0} q_i^*
q_i\Bigr\|.$$
It is clear that
$$\Bigl\| \sum_{i \in I} p_i \ot q_i \Bigr\|_h^2 \leq
\Bigl\| \sum_{i \in I} p_i p_i^* \Bigr\| \, \Bigl\| \sum_{i \in I} q_i^*
q_i \Bigr\|.$$
It is also possible to prove a converse of this statement: every
element of $A \ot_h A$ can be written as such an infinite sum. We
do not go into this because we will not need this result.

Quite often we will extend maps defined on the algebraic tensor
product $A \od A$ to $A \oh A$ by showing their continuity. We
will use the following obvious lemma (see also lemma~\ref{41b}).
\begin{lemma} \label{11b}
Let $B$ be a C$^*$-algebra and $p_i,q_i \in B$ for $i=1,\ldots,n$.
Then
$$\Bigl\| \sum_{i=1}^n p_i q_i \Bigr\|^2 \leq \Bigl\|\sum_{i=1}^n p_i p_i^* \Bigr\| \,
\Bigl\| \sum_{i=1}^n q_i^* q_i \Bigr\|.$$
\end{lemma}
Now we prove the following crucial result.
\begin{prop} \label{12}
The inclusion of $A \od A$ in $A \ot A$ is a contraction for the
Haagerup norm on $A \od A$ and can be extended uniquely to a contraction from $A
\oh A$ to $A \ot A$. This extension is injective.
\end{prop}
\begin{proof}
When $p_i,q_i \in A$ we have by the previous lemma
$$\Bigl\|\sum_{i=1}^n p_i \ot q_i \Bigr\|^2 = \Bigl\|\sum_{i =1}^n (p_i \ot 1) (1 \ot
q_i) \Bigr\|^2 \leq \Bigl\|\sum_{i=1}^n p_i p_i^* \Bigr\| \,
\Bigl\| \sum_{i=1}^n q_i^* q_i \Bigr\|.$$
So we can extend the inclusion of $A \od A$ in $A \ot A$ to a
contraction $j$ from $A \oh A$ to $A \ot A$. Now we have to prove the
injectivity of $j$. Suppose $j(x)=0$ for $x \in A \oh A$. Then we
have $(\om \ot \mu)(x)=(\om \ot \mu)j(x)=0$ for all $\om, \mu \in A^*$. Because
the inclusion of $A \oh A$ in $(A^* \oh A^*)^*$ is completely
isometric (see e.g.~\cite{E-R}), we may conclude that $x = 0$.
\end{proof}
From now on we will identify $A \oh A$ with a dense subspace of $A
\ot A$.

From lemma~\ref{11b} follows immediately the following result.
\begin{prop} \label{13}
The multiplication map $m : p \otimes q \mapsto
pq$ can be uniquely extended to a contraction from the Haagerup tensor product $A \otimes_h A$
to $A$.
\end{prop}
We know that the slice maps $\iota \otimes \omega$ and $\omega \otimes \iota$ are
continuous from $A \otimes A$ to $A$ and that the norms of these slice maps are
majorized by $\| \omega \|$.  This is of course also true for the restriction
of these slice maps to $A \otimes_h A$. More generally, let $A,B$ and $C$ be C$^*$-algebras.
When $T : A \rightarrow B$ is a so-called completely bounded map,
then the map $T \ot \io$ from $A \od C$ to $B \od C$ is bounded
for the Haagerup norm. Its norm is majorized by $\|T\|_{cb}$, the
completely bounded norm of $T$. For this we refer to the discussion in the beginning of section 4.
Examples of such completely bounded maps are slice maps $\om \ot
\io$ from the C$^*$-algebra $B \ot A$ to the C$^*$-algebra $A$, or $^*$-homomorphisms between
C$^*$-algebras.

As a special case of this (but that is of course more easy to see) we have that
the maps
$$x \mapsto x(1 \otimes a) \quad  {\rm and } \quad x \mapsto (a \otimes 1)x$$
are bounded from $A \otimes_h A$ to itself with norms majorized by $\| a \|$.
This, together with the following lemma will be the starting point for the
introduction of a larger space, the so-called multiplier Haagerup tensor
product.

\begin{lemma} \label{14}
If $(e_\lambda)$ is a bounded approximate identity in $A$,
then $x(1 \otimes e_\lambda) \rightarrow x$ and $(e_\lambda \otimes 1)x
\rightarrow x$ in $A \otimes_h A$ for all $x \in A \otimes_h A$.
\end{lemma}
\begin{proof}
It is sufficient to do this for elements in $A \odot A$.  If $x
= \sum_{i=1}^n p_i \otimes q_i$ then $x -x(1 \otimes e_\lambda) = \sum_{i=1}^n p_i \otimes
q_i(1-e_\lambda)$ and so
$$\| x - x(1 \otimes e_\lambda) \|_h \leq \Bigl\| \sum_{i=1}^n p_i p_i^\ast \Bigr\|^{1/2} \Bigl\|
\sum_{i=1}^n (1-e_\lambda^\ast)q^\ast_i q_i(1 - e_\lambda) \Bigr\|^{1/2}$$
and this last expression will tend to $0$.
\end{proof}
We will now introduce the multiplier Haagerup tensor product $A \otimes_{mh}A$.
The way we will do this here is different from~\cite{V}.  In principle we get a
larger space, although we are not completely sure about that.  In any case, this
new approach seems to be somehow more natural.  We define $A \otimes_{mh} A$
in a way which is very similar to the definition of the multiplier algebra $M(A
\otimes A)$ of $A \otimes A$.  We will see later in the next section why this
is natural.  In fact, this should also become clear from certain properties
that we will prove already here.

\begin{defin} \label{15}
Denote by $A \otimes_{mh}A$ the space of pairs of maps
$(\rho_1,\rho_2)$ from $A$ to $A \otimes_h A$ satisfying $(a \otimes
1)\rho_1(b) = \rho_2(a)(1 \otimes b)$ for all $a,b \in A$.
\end{defin}

We clearly have $A \otimes_h A$ sitting in $A \otimes_{mh} A$ by
putting
$\rho_1(a) = x(1 \otimes a)$ and $\rho_2 (a) = (a \otimes 1)x$ where $x \in A
\otimes_h A$ and $a \in A$.  We also have an imbedding as $x(1 \otimes a) = 0$
for all $a$ will imply $x=0$ by the previous lemma.

We will  now show that $A \otimes_{mh} A$ carries a norm such that this space is
complete and contains $A \otimes_h A$ isometrically.

\begin{prop} \label{16}
Let $x = (\rho_1,\rho_2) \in A \otimes_{mh}A$.  Then
$\rho_1$ and $\rho_2$ are bounded and $\| \rho_1 \| = \| \rho_2 \|$.
Further $\| x \| := \| \rho_1 \| = \| \rho_2 \|$ defines a norm on $A \otimes_{mh} A$
for which this space is complete.  Moreover, the imbedding of $A \otimes_h A$
in $A \otimes_{mh} A$ described above is isometric.
\end{prop}

\begin{proof}
We will first show that the maps $\rho_1$ and $\rho_2$ are
closed.  Take $a_n \in A$ and assume that $a_n \rightarrow a$ in $A$ and
$\rho_1 (a_n) \rightarrow x$ in $A \otimes_h A$.  Then for any $b \in A$ we
have
\begin{align*}
(b \otimes 1)x & = \lim_n (b \otimes 1)\rho_1(a_n) \\
& = \lim_n \rho_2 (b) (1 \otimes a_n) \\
& = \rho_2 (b)(1 \otimes a)\\
& = (b \otimes 1)\rho_1(a).
\end{align*}
Then $x = \rho_1(a)$.  By the closed graph theorem we have that $\rho_1$ is
bounded.  Similarly $\rho_2$ is bounded.  Now, because of
lemma~\ref{14} we have
\begin{align*}
\| \rho_1 \| & = \sup \{ \| \rho_1(a) \|_h \mid \| a \| \leq 1 \} \\
& = \sup \{ \| (b \otimes 1) \rho_1 (a) \|_h \mid \| a \| \leq 1, \| b \| \leq 1
\} \\
& = \sup \{ \| \rho_2 (b) (1 \otimes a) \|_h \mid \| a \| \leq 1, \| b \| \leq 1
\} \\
& = \| \rho_2 \| .
\end{align*}
As for $x \in A \otimes_h A$ we have $\| x \|_h = \sup \{  \| x(1 \otimes a)
\|_h \mid \| a \| \leq 1 \}$ we get that the imbedding of $A \otimes_h A$ in $A
\otimes_{mh} A$ is isometric.

Finally we show that $A \otimes_{mh} A$ is complete.  Take a Cauchy sequence
$(x_n)$ in $A \otimes_{mh} A$.  Write $x_n = (\rho_{1n}, \rho_{2n})$.  Then
$(\rho_{1n}(a))_n$ and $(\rho_{2n} (a))_n$ are Cauchy sequences in $A
\otimes_h
A$.  They have limits and we get maps $(\rho_1,\rho_2)$ defined by the
pointwise limit.  Obviously they still satisfy $(a \otimes 1)\rho_1(b)  =
\rho_2(a)(1 \otimes b)$ for all $a,b \in A$.  By a standard technique we get
that $\rho_{1n} \rightarrow \rho_1$ and $\rho_{2n} \rightarrow \rho_2$ also in
norm.  This completes the proof.
\end{proof}

Further we will write $x(1 \otimes a)$ for $\rho_1(a)$ and $(a \otimes 1)x$ for
$\rho_2(a)$, when $x = (\rho_1, \rho_2)$ for any $x \in A \otimes_{mh} A$.
This notation is compatible with the relation between $\rho_1$ and  $\rho_2$ (now
expressed as associativity), with the imbedding of $A \otimes_h A$ in $A
\otimes_{mh} A$ and with the norm properties.  Moreover, we can justify this
notation even further because we have the following result.

\begin{prop} \label{17}
The inclusion $j : A \oh A \rightarrow A \ot A$ can be uniquely
extended to a contraction from $A \omh A$ to $M(A \ot A)$, still
denoted by $j$, such that
$$j(x) (1 \ot a) = j(x(1 \ot a)) \quad \text{and} \quad (a \ot 1)
j(x) = j((a \ot 1)x)$$
for all $x \in A \omh A$ and $a \in A$. Also this extension is
injective.
\end{prop}
\begin{proof}
Take $x \in A \otimes_{mh} A$.  As $x (1 \otimes a) \in A
\otimes_h A \subseteq A \otimes A$ we can clearly define a left multiplier from
$A \odot A$ to $A \otimes A$ by $x(b \otimes a) = x(1 \otimes a)(b \otimes
1)$.  We will show that this is bounded and therefore extends to a left
multiplier from $A \otimes A$ to $A \otimes A$.  To see this, take a bounded
approximate unit $(e_\lambda)$ of $A$ and observe that
$$(e_\lambda \otimes 1) \Bigl(x \Bigl(\sum_{i=1}^n b_i \otimes a_i \Bigr)\Bigr) =
((e_\lambda \otimes 1)x)
\Bigl(\sum_{i=1}^n b_i  \otimes a_i\Bigr).$$
Now $\| (e_\lambda \otimes 1)x \| \leq \| (e_\lambda \otimes 1)x \|_h \leq \| x
\|_h$ and so
$$\Bigl\| x \Bigl( \sum_{i=1}^n b_i \otimes a_i \Bigr) \Bigr\| \leq \| x \|_h \,
\Bigl\| \sum_{i=1}^n b_i \otimes a_i \Bigr\|.$$
Similarly, $x$ defines a bounded right multiplier of $A \otimes A$.  And
obviously $((a \otimes b)x)(c \otimes d) = (a \otimes b)(x(c \otimes d))$ so
that indeed $x \in M(A \otimes A)$.
\end{proof}
From now on we will identify $A \omh A$ with a subspace of $M(A
\ot A)$.
So we have $A \otimes_{mh} A$ in $M(A \otimes A)$ and now $x(1 \otimes a)$
makes sense in $M(A \otimes A)$.  In fact we get the following
characterization.
\begin{prop} \label{18}
If $x \in M(A \otimes A)$, then $x \in A \otimes_{mh}
A$ if and only if $x(1 \otimes a) \in A \otimes_h A$ and $(a \otimes 1)x \in A
\otimes_h A$ for all $a \in A$.
\end{prop}
It is clear that we could also have started with this property to define $A
\otimes_{mh}A$ from the beginning as sitting in $M(A \otimes A)$.

Now, very often,  we will need to extend maps from $A \otimes_h A$ to $A
\otimes_{mh} A$.  The first example we have in mind is multiplication :
\begin{prop} \label{19}
The multiplication map $m$ can be extended uniquely to a contraction from
$A \otimes_{mh} A$ to $M(A)$ satisfying
\begin{align*}
m(x)a &= m(x(1 \otimes a)) \\
a m(x) &= m((a \otimes 1)x)
\end{align*}
when $x \in A \otimes_{mh} A$ and $a \in A$.
\end{prop}
\begin{proof}
Now, if $x \in A \otimes_{mh} A$ we can define two maps from $A$
to $A$ by
\begin{align*}
a & \mapsto m(x(1 \otimes a)) \\
a & \mapsto m((a \otimes 1) x)
\end{align*}
where we use now $m$ for the multiplication from $A \otimes_h A$ to $A$.  Because
\begin{align*}
a m(x(1 \otimes b)) & = m ((a \otimes 1)x(1 \otimes b)) \\
& = m((a \otimes 1)x)b
\end{align*}
this pair of maps gives an element in $M(A)$.
\end{proof}
This is the kind of technique we will be using, but it is good to know what is
really behind.  This is what we will explain now.
\begin{defin} \label{110}
By the strict topology on $A \otimes_{mh} A$ we mean
the locally convex topology given by the semi-norms $x \mapsto \| x(1
\otimes a) \|_h$ and $x \mapsto \| (a \otimes 1) x \|_h$ for $a \in A$.
\end{defin}
We have the following density result.
\begin{prop} \label{111}
When $x \in A \otimes_{mh} A$ and $(e_\lambda)$ is a
bounded approximate identity in $A$, then $x(1 \otimes e_\lambda) \rightarrow
x$ and $(e_\lambda \otimes 1)x \rightarrow x$ strictly.
\end{prop}
\begin{proof}
Take $a \in A$ and observe that $(a \otimes 1)(x(1 \otimes e_\lambda)) =
((a \otimes 1)x)(1 \otimes e_\lambda) \rightarrow (a \otimes 1)x$ in $A
\otimes_h A$ (lemma~\ref{14}).  On the other hand also $x(1 \otimes e_\lambda)(1 \otimes
a) = x(1 \otimes e_\lambda a) \rightarrow x(1 \otimes a)$ as $e_\lambda a \rightarrow a$ in
norm in $A$.
\end{proof}
It follows that the unit ball of $A \otimes_h A$ is strictly dense in
the unit ball of $A \otimes_{mh} A$.  Therefore to extend maps one
can show that they are continuous for the strict topology on bounded sets.
Such is e.g.\ the case for the multiplication map.  In fact, $m$ as defined
in~\ref{19} from $A \otimes_{mh}A$ to $M(A)$ is continuous when both spaces are taken
with the strict topology. From the proof of proposition~\ref{16}
also follows immediately that the inclusion of $A \omh A$ in $M(A
\ot A)$ is strictly continuous on bounded sets.

Also slice maps are naturally defined on $A \otimes_{mh} A$:
\begin{prop}
Let $\om \in A^*$. The maps $\om \ot \io$ and $\io \ot \om$ from
$A \oh A$ to $A$ extend uniquely to strict-norm continuous maps
from $A \omh A$ to $A$ satisfying
$$a (\om \ot \io)(x) = (\om \ot \io)((1 \ot a)x) \quad \text{and}
\quad (\om \ot \io)(x) a = (\om \ot \io)(x(1 \ot a))$$
and satisfying analogous formulas for $\io \ot \om$. The norms of these maps
are majorized by $\|\om\|$.
\end{prop}
Remark that we mean by strict-norm continuity from $A \omh A$ to
$A$ continuity when $A \omh A$ is equipped with the strict
topology and $A$ with the norm topology.
\begin{proof}
This is quite easy to see. We look at the case of $\om \ot \io$.
Then the formulas in the proposition define an element in $M(A)$.
Now write $\om = \mu b$ with $\mu \in A^*$ and $b \in A$ and
observe that $(\om \ot \io)(x)= (\mu \ot \io)((b \ot 1)x) \in A$. This
gives the result.
\end{proof}
There is one more result on this space we will need.
\begin{prop} \label{112}
The map $\Gamma : p \otimes q \rightarrow q^\ast
\otimes p^\ast$ extends to an isometry of $A \otimes_{mh} A$ satisfying
\begin{equation*}
\Gamma(x(1 \otimes a))  = (a^\ast \otimes 1) \Gamma (x) \quad
\quad \Gamma((a \otimes 1)x)  = \Gamma (x)(1 \otimes a^\ast)
\end{equation*}
for $x \in A \omh A$ and $a \in A$.
The map $\Gamma$ is a homeomorphism for the strict topology on $A
\omh A$.
\end{prop}
\begin{proof}
First we need to verify that $\Gamma$ leaves the Haagerup norm
invariant.  This is more or less straightforward.  Then $\Gamma$ is extended
to $A \otimes_h A$.  Finally the formulas in the proposition will yield the
extension of $\Gamma$ to $A \otimes_{mh} A$ and also show the
strict continuity.
\end{proof}
The following property of $A \otimes_{mh} A$ will not really be needed.  One
can show that $A \otimes_{mh} A$ sits in the extended Haagerup tensor product
$A \otimes_{eh} A$ as defined in~\cite{EKR}.  It follows that elements in $A
\otimes_{mh} A$ can be written in the form $\sum p_i \otimes q_i$ with $\sum
p_i p_i^\ast$ and $\sum q^\ast_i q_i$ weakly convergent to some elements in the
double dual $A^{**}$.  The sum $\sum p_i \otimes q_i$ is now convergent in the
von Neumann algebra tensor product $A^{**} \bar{\otimes} A^{**}$.

This takes us to the following natural subspace of $A \otimes_{mh} A$.
\begin{prop} \label{113}
Let $(p_i),(q_i)$ be families of elements of $A$ such that the
$\sum_{i \in I} p_i p_i^\ast$ and $\sum_{i \in I} q_i^\ast q_i$ converge strictly
and remain bounded in $M(A)$, then the net $\sum_{i \in I}
p_i \otimes q_i$ converges strictly to an element in $A \otimes_{mh} A$.
\end{prop}
\begin{proof}
Let $a \in A$. Then the net $\sum_{i \in I} a p_i p_i^* a^*$
is norm convergent, and the net $\sum_{i \in I} q_i^* q_i$
is bounded. So, just as in the remark after definition~\ref{11}
you get that the net
$$\sum_{i \in I} a p_i \ot q_i$$
is norm convergent in $A \oh A$. We denote the limit by
$\rho_2(a)$. Analogously we define
$$\rho_1(a) = \sum_{i \in I} p_i \ot q_i a.$$
Then it is clear that $x=(\rho_1,\rho_2)$ belongs to $A \omh A$, and
by definition the net $\sum_{i \in I} p_i \ot q_i$ converges strictly to
$x$.
\end{proof}
In~\cite{V}, the present theory has been developed using this subspace of $A
\otimes_{mh} A$. We do not know whether this subspace is the whole
of $A \otimes_{mh} A$ or not.

\section{Hopf C$^*$-algebras}
As before, let $A$ be a C$^*$-algebra.  We consider the minimal C$^*$-tensor product
$A \otimes A$ of $A$ by itself.  Let $\phi$ be a non-degenerate
$^\ast$-homomorphism of $A$ into the multiplier algebra $M(A \otimes A)$ of $A
\otimes A$.  Denote by $\iota$ the identity map of $A$ and consider the maps
$\iota \otimes \phi$ and $\phi \otimes \iota$
from $A \otimes A$ to $M(A \otimes A \otimes A)$.  They are again
non-degenerate and therefore have unique extensions to $^\ast$-homomorphisms
from $M(A \otimes A)$ to $M(A \otimes A \otimes A)$.  These extensions will still be
denoted by $\iota \otimes \phi$ and $\phi \otimes \iota$.

Then we are ready for the following definition.
\begin{defin} \label{21}
A comultiplication on a C$^*$-algebra $A$ is a
non-degenerate $^\ast$-homomorphism $\phi : A \rightarrow M(A \otimes A)$ such
that $(\phi \otimes \iota) \phi = (\iota \otimes \phi)\phi$.
\end{defin}
The condition of
coassociativity, namely $(\phi \otimes \iota)\phi = (\iota \otimes
\phi)\phi$,
is given a meaning by extending $\phi \otimes \iota$ and $\iota \otimes \phi$
to $M(A \otimes A)$ as above.  For this it is necessary to have $\phi$
non-degenerate.

The above conditions are satisfied in the following basic example.
\begin{prop} \label{22}
Let $G$ be a locally compact semi-group.  Let $A$ be the
C$^*$-algebra $C_0(G)$ of continuous complex functions on $G$ tending to $0$
at infinity.  Then $\phi$ defined by $\phi(f)(p,q) = f(pq)$ when $f \in
C_0(G)$ and $p,q \in G$, is a comultiplication on $A$.
\end{prop}
\begin{proof}
First observe that $A \otimes A$ is naturally identified with
$C_0(G \times G)$ and that $M(A \otimes A)$ is hence identified with $C_b(G
\times G)$, the C$^*$-algebra of bounded continuous complex functions on $G
\times G$.  Therefore, $\phi$ as defined in the formulation of the proposition
is a map from $A$ to $M(A \otimes A)$.  It is clearly a $^\ast$-homomorphism
and it is not hard to show that $\phi$ is non-degenerate.  The
coassociativity comes from the associativity of the product in $G$.
\end{proof}
It follows that we can think of a pair of a C$^*$-algebra $A$ and a
comultiplication $\phi$ on $A$ as a locally compact quantum semi-group.  Now,
as already explained in the introduction, it is not at all clear how to impose
extra conditions in order to get a locally compact quantum group.  One cannot
simply require the existence of a counit and an antipode.

The most ambitious project would be to look for extra conditions that are easy
to verify in special cases (such as the locally compact groups) and nice
non-trivial examples and such that a theory can be developed along the lines of
locally compact groups.  In such a frame, the existence of (the equivalent of)
the Haar measure would be one of the main results.  It seems however that such
an objective is still far out of reach.

All the known existence proofs for the Haar measure in the classical case seem
to fail in the non-commutative setting (except for special cases like the
compact and the discrete quantum groups).  Therefore, it is quite acceptable to
develop a theory where the existence of the Haar measure is assumed.  After
all, just as in the classical case, also in the quantum case the Haar measures
are simply there. From this point of view, the uniqueness property is at least
as important.

Such a theory has now, very recently, been developed and the results are very
satisfactory.  In~\cite{KV1} and \cite{KV2}, a pair $(A, \phi)$ of a C$^*$-algebra with a
comultiplication is called a locally compact quantum group if the spaces
$\phi(A)(1 \otimes A)$ and $(A \otimes 1)\phi(A)$ are dense in $A \otimes A$
and if there exist 'good' left and right Haar measures.  Remark that these
density requirements correspond in the classical case (with $A$ abelian) with
the cancellation law. Also observe that in case $A$ has an identity, the
density conditions are sufficient to prove the existence of a Haar measure.
This gives the theory of compact quantum groups.  But of course, this will not
be enough in general (when $A$ has no identity).

Nevertheless, the question remains what is a good notion for a counit and an
antipode in this context.  In fact, these notions do play an important role in
the theory of locally compact quantum groups of~\cite{KV2}.  Moreover, it is
worthwhile to develop such notions independently of the existence of Haar
measures.  This is what we will do further in the paper.  And to distinguish
from the theory of locally compact quantum groups with Haar measure, we will
call the objects we get here Hopf C$^*$-algebras.  This name is appropriate
for a C$^*$-algebra with a comultiplication, a counit and an
antipode. Cfr.\ the discussion concerning this terminology in the
introduction.

To see what extra conditions might be natural, let us think of the theory of
Hopf algebras.  If $(H,\Delta,\varepsilon,S)$ is a Hopf algebra, then the maps
$T_1, T_2$, defined on $H \otimes H$ by
\begin{align*}
T_1 (a \otimes b) & = \Delta(a)(1 \otimes b) \\
T_2 (a \otimes b) & = (a \otimes 1)\Delta(b)
\end{align*}
are bijective.  The inverses are given in terms of the antipode by
\begin{align*}
T^{-1}_1(a \otimes b) & = (( \iota \otimes S)\Delta(a))(1 \otimes b) \\
T^{-1}_2(a \otimes b) & = (a \otimes 1)((S \otimes 1)\Delta(b)).
\end{align*}
In multiplier Hopf algebra theory, the bijectivity of these maps was the
starting point from which the counit and the antipode where constructed (see~\cite{VD2}).

If we want to develop a C$^*$-algebra theory in this direction, it is more
or less obvious
what we have to do : we need to require that the maps $T_1$ and $T_2$ above
have dense range in $A \otimes A$ and that proper extensions are (still)
injective.

The abelian case here is very misleading.  If $A = C_0(G)$ with $G$ a locally
compact group, then the maps $T_1$ and $T_2$ extend to $^\ast$-isomorphisms
of $A \otimes A$ given by $(T_1 f)(p,q) = f(pq,q)$ and $(T_2f)(p,q) = f(p,pq)$
when $f \in C_0(G \times G)$ and $p,q \in G$.  This is not at all the case in
general where e.g.\ $T_1$ and $T_2$ will not even be bounded (the reason being
that multiplication is no longer bounded on $A \otimes A$ with the
C$^*$-norm).

Fortunately, as we have seen in the previous section, multiplication is
continuous for the Haagerup norm.  And it is precisely that property of the
Haagerup norm that allows the following result which is the starting point in
the development of our notion of a Hopf C$^*$-algebra.
\begin{prop} \label{23}
Let $(A,\phi)$ be a pair of a C$^*$-algebra and a comultiplication
such that $\phi(A)(1 \ot A)$ and $\phi(A)(A \ot 1)$ are subspaces
of $A \ot A$.
Then the maps $T_1$ and $T_2$ defined from $A \odot A$
to $A \otimes A$ by
\begin{equation*}
T_1(a \otimes b)  = \phi(a)(1 \otimes b) \quad \quad
T_2(a \otimes b)  = (a \otimes 1) \phi (b)
\end{equation*}
are contractive from $A \odot A$ with the Haagerup norm to $A
\otimes A$ with the C$^*$-norm. We use the same symbols $T_1$ and
$T_2$ to denote the extensions to $A \oh A$.
\end{prop}
\begin{proof}
Take $x \in A \odot A$ and write $x  = \sum_{i=1}^n p_i \otimes q_i$.
By using lemma~\ref{11b} we get
\begin{align*}
\|T_1(x)\| &= \Bigl \| \sum_{i=1}^n \phi (p_i)(1 \otimes q_i) \Bigr\| \\
& \leq \Bigl\| \sum_{i=1}^n \phi(p_i p_i^*) \Bigr\|^{1/2} \, \Bigl\| \sum_{i=1}^n
q_i^* q_i \Bigr\|^{1/2} \\
& \leq \Bigl\|  \sum_{i=1}^n p_i p_i^* \Bigr\|^{1/2} \, \Bigl\| \sum_{i=1}^n
q_i^* q_i \Bigr\|^{1/2}.
\end{align*}
So we can conclude that $\|T_1(x)\| \leq \|x\|_h$. Analogously
$\|T_2(x)\| \leq \|x\|_h$.
\end{proof}
We see that we use essentially the same argument as in the proof of
proposition~\ref{13}
where we show that multiplication is contractive.  This is no surprise
as we can write $T_1 = (\iota \otimes m)(\phi \otimes \iota)$ and $T_2 = (m
\otimes \iota)(\iota \otimes \phi)$.  Remark that the order of $a$ and $b$ is
important (the Haagerup norm is not symmetric, i.e.\ not invariant under the
flip).  This means that a map like $a \otimes b \mapsto (1 \otimes
b)\phi(a)$ will in general not be continuous for the Haagerup norm, but a map
like $a \otimes b \mapsto (1 \otimes a)\phi (b)$ will be.  We will use this
last map also in section 3.

Then we are ready for the main definition.
\begin{defin} \label{24}
Let $(A, \phi)$ be a pair of a C$^*$-algebra $A$
with a comultiplication $\phi$ such that $\phi(A)(1 \ot A)$ and $\phi(A)(A \ot 1)$ are subspaces
of $A \ot A$.  We call this pair a Hopf C$^*$-algebra if
the maps $T_1$ and $T_2$ as defined above, are injective on $A \otimes_h A$.
\end{defin}
As we have said before, the injectivity of $T_1$ and $T_2$ on $A \otimes_h A$
is obvious in the case of $A = C_0(G)$ with $G$ a locally compact group (where we have much stronger
properties).  But what is rather striking is that also conversely, all abelian
Hopf C$^*$-algebras are of this form. We even have the following
stronger result.

\begin{prop}
Let $(A,\phi)$ be a pair of an abelian C$^*$-algebra $A$ with a
comultiplication $\phi$ such that $\phi(A)(1 \ot A)$ and $\phi(A)(A \ot 1)$ are subspaces
of $A \ot A$. Assume further that
\begin{align*}
\phi(a)(1 \ot b) = 0 &\quad\text{implies}\quad a \ot b = 0 \\
\text{and}\quad \phi(a)(b \ot 1) = 0 &\quad\text{implies}\quad a \ot b = 0
\end{align*}
for all $a,b \in A$. Then there exists a locally compact group $G$
such that $A \cong C_0(G)$ and $\phi$ is given by
$\phi(f)(p,q)=f(pq)$ for all $p,q \in C_0(G)$ under this
isomorphism.
\end{prop}
\begin{proof}
It is clear by the Gelfand-Naimark theorem that there exists a
locally compact space $G$ and a continuous associative binary operation on $G$
such that $A \cong C_0(G)$ and such that $\phi$ is given by
the formula above under this isomorphism. We will prove now that
$G$ is a group and that the inverse is continuous.

First choose a non-empty open subset $V$ of $G$. We claim that
the set $\{rs \mid r \in G, s \in V\}$ is dense in $G$. Indeed,
suppose that it is not dense. Then there is an element $f \in
C_0(G)$ such that $f \neq 0$ but $f(rs)=0$ for all $r \in G$ and
$s \in V$. Choose an element $g \in C_0(G)$ such that $g \neq 0$
and such that $g$ has support in $V$. Then $f(rs)g(s)=0$ for all
$r,s \in G$. Hence $\phi(f)(1 \ot g)=0$ and by assumption we may
conclude that $f=0$ or $g = 0$. So we get a contradiction and this
gives the claim.

Now let $s,t \in G$. For any pair $\alpha=(V,W)$ of open neighborhoods of
$s$ and $t$ respectively, we have, by the property above,  points
$r_\alpha \in G$, $s_\alpha \in V$ and $t_\alpha \in W$ such that
$r_\alpha s_\alpha = t_\alpha$. These pairs of neighborhoods give
an index set $I$ when ordered by opposite inclusions. By
definition the net $(s_\alpha)$ converges to $s$ and $(t_\alpha)$
to $t$.

For every $f,g \in C_0(G)$ we have
$$(\phi(f)(1 \ot g))(r_\alpha,s_\alpha) = f(r_\alpha s_\alpha)
g(s_\alpha) = f(t_\alpha) g(s_\alpha)$$
and this converges to $f(t)g(s)$. Now choose $f$ and $g$ such that
$f(t)g(s) \neq 0$. Because $\phi(f)(1 \ot g) \in C_0(G \times G)$
by assumption, we can take compact sets $K_1$ and $K_2$ such that
$|\phi(f)(1 \ot g)|$ is smaller then $\frac{1}{2} |f(t)g(s)|$
outside $K_1 \times K_2$. Then, for $\alpha$ large enough, we will
have $|(\phi(f)(1 \ot g))(r_\alpha,s_\alpha)|  \geq \frac{1}{2}
|f(t)g(s)|$ and so $r_\alpha \in K_1$. Then we can find a subnet
of $(r_\alpha)$ that converges to a point $r$. Hence we find $r
\in G$ such that $rs=t$.

Because $s,t \in G$ were chosen arbitrarily we get $G s = G$ for
all $s \in G$. By symmetry we also get $s G = G$ for all $s \in
G$. Then it is an easy exercise to prove that $G$ is a group. In
particular we have the cancellation property. Therefore the linear
span of $\{\phi(f)(1 \ot g) \mid f,g \in C_0(G) \}$ is a $^*$-subalgebra
of $C_0(G \times G)$ that separates points of $G \times G$ and
hence is dense in $C_0(G \times G)$. Now for any $f,g \in C_0(G)$
the map
$$s \mapsto (\phi(f)(1 \ot g))(s^{-1}, s) = f(s^{-1}s) g(s) = f(e)
g(s)$$
is continuous. By density we get that the map $s \mapsto f(s^{-1})
g(s)$ is continuous for all $f,g \in C_0(G)$. From this it follows
that $s \mapsto s^{-1}$ is continuous and that concludes the
proof.
\end{proof}

This proposition is one reason for not assuming the density of the ranges of $T_1$ and
$T_2$ in the definition of a Hopf C$^*$-algebra,
which is a little surprising after the discussion we have
had.  Another reason is that on the one hand we do not need these density
conditions for the development of our theory,
while on the other hand we will formulate another
density condition that will be stronger.  We will come back to this in
section 3.

A problem with the notion of a Hopf C$^*$-algebra is that injectivity on this completed space $A
\otimes_h A$ may sometimes be hard to verify.  We will however give
techniques to do this
in section 5, where we treat examples and special cases.  There we
will see that we do get Hopf C$^*$-algebras in most interesting cases.

Now we need also the extension of $T_1$ and $T_2$ to $A \otimes_{mh} A$.  The
proof that these extensions exist, is essentially similar to the fact that also
multiplication can be extended to $A \otimes_{mh} A$.

\begin{prop} \label{25}
The maps $T_1$ and $T_2$ can be uniquely extended to
contractions from $A \otimes_{mh} A$ to $M(A \otimes A)$
satisfying
\begin{align*}
T_1(x)(1 \otimes a) & = T_1 (x(1 \otimes a)) && \text{and} & \phi(a)T_1(x) &= T_1((a \ot 1)x) \\
(a \otimes 1)T_2(x) & = T_2((a \otimes 1)x) && \text{and}
& T_2(x) \phi(a) &= T_2(x(1 \ot a))
\end{align*}
whenever $x \in A \otimes_{mh} A$ and $a \in A$.  These extensions are still
injective and they are strictly continuous on bounded sets.
\end{prop}
\begin{proof}
The injectivity of the extensions is an immediate consequence of
these defining formulas.  If e.g.\ $T_1(x) = 0$ for $x \in A \otimes_{mh}A$,
then $x(1 \otimes a) = 0$ for all $a \in A$ by the injectivity of $T_1$ on $A
\otimes_h A$ and hence $x = 0$.  So it essentially remains to show that the
above formulas yield maps into $M(A \otimes A)$.

So, let $x \in A \otimes_{mh} A$.  We define a left multiplier from $A \odot
A$ to $A \ot A$
by $T_1 (x) (b \otimes a) = T_1 (x(1 \otimes a))(b \otimes 1)$.  When
$(e_\lambda)$ is an approximate identity in $A$, we have
\begin{align*}
\phi(e_\lambda)T_1(x) \Bigl(\sum_{i=1}^n b_i \otimes a_i \Bigr) & = \phi (e_\lambda) \Bigl(\sum_{i=1}^n T_1 (x (1
\otimes a_i))(b_i \otimes 1)\Bigr) \\
& = \sum_{i=1}^n T_1 ((e_\lambda \otimes 1)x(1 \otimes a_i))(b_i \otimes 1) \\
& = T_1((e_\lambda \otimes 1)x)\Bigl(\sum_{i=1}^n b_i \otimes a_i\Bigr).
\end{align*}
This implies that
\begin{align*}
\Bigl\| \phi(e_\lambda)T_1(x)\Bigl(\sum_{i=1}^n b_i \otimes a_i\Bigr) \Bigr\|  & \leq \| (e_\lambda \otimes 1)x
\|_h \, \Bigl\| \sum_{i=1}^n b_i \otimes a_i \Bigr\| \\
& \leq \| x \|_h \, \Bigl\| \sum_{i=1}^n b_i \otimes a_i \Bigr\|.
\end{align*}
Taking the limit we get
$$\Bigl\| T_1(x)\Bigl(\sum_{i=1}^n b_i \otimes a_i \Bigr) \Bigr\| \leq \| x \|_h \, \Bigl\| \sum_{i=1}^n b_i \otimes a_i
\Bigr \|.$$
This shows that we have a left multiplier of $A \otimes A$ with norm less than
$\| x \|_h$.  Similarly, we define $T_1 (x)$ as a right multiplier by
$$\Bigl(\sum_{i=1}^n (a_i \otimes b_i)\phi(c_i)\Bigr) T_1(x) = \sum_{i=1}^n
(a_i \otimes b_i)T_1((c_i \otimes 1)x).$$
A similar argument as above, by looking at $x(1 \otimes e_\lambda)$, together with the
non-degeneracy of $\phi$, will give that $T_1 (x)$ is also a bounded right
multiplier of $A \otimes A$.

Hence $T_1(x) \in M(A \otimes A)$ and $\| T_1(x) \| \leq \| x
\|_h$. From the formulas above it follows immediately that $T_1$
extended to $A \omh A$ will be strictly continuous on bounded
sets. We can give a similar argument for $T_2$.
\end{proof}

\section{The left counit and antipode}
Now let $A$ be a C$^*$-algebra with a comultiplication $\phi$ satisfying the
assumptions of the previous section, making $(A,\phi)$ a Hopf C$^*$-algebra.  So, the maps $T_1$ and $T_2$ defined
from
$A \odot A$ to $A \otimes A$ by $T_1 (a \otimes b) = \phi(a)(1 \otimes b)$ and $T_2(a \otimes
b) = (a \otimes 1)\phi(b)$ have ranges in $A \otimes A$ and their
continuous
extensions to $A \otimes_h A$ are injective.  Again as in the previous section,
we consider their extensions to $A \otimes_{mh} A$, still using the same
symbols.  These extensions are still injective.  Now the range is in $M(A
\otimes A)$.  Recall that the maps $T_1,T_2$ from $A \otimes_{mh} A$ to $M(A
\otimes A)$ are strictly continuous on bounded sets.

We will now first define a candidate for the domain of the counit and the
antipode.
\begin{nota} \label{31}
Let $A_0$ be the space of elements $a \in A$ such that
there is an element $x \in A \otimes_{mh} A$ satisfying $a \otimes 1 = T_1x$
\end{nota}
In case $A = C_0(G)$ with $G$ a locally compact group, it is possible to show
that $A_0$ is dense.  Also for compact and discrete and more generally locally compact quantum groups this is
the case. We will prove these things in section~5.
However, without further assumptions, it seems that nothing
can be said about the size of $A_0$. Nevertheless it is
immediately clear that the density of $A_0$ implies that
$\phi(A)(1 \ot A)$ is dense in $A \ot A$, so that we will
have to assume at least that to get $A_0$ dense in $A$.

We will not get more elements if we would allow $a \in M(A)$. Indeed, if $a \in
M(A)$ and $x \in A \otimes_{mh} A$ satisfy $a \otimes 1 = T_1 x$, then $a
\otimes b = (T_1 x)(1 \otimes b) = T_1(x(1 \otimes b))$ which belongs to $A
\otimes A$ for all $b \in A$ and this will imply $a \in A$.

It is obvious that $A_0$ is a subspace of $A$.  But we have more :
\begin{prop} \label{32}
$A_0$ is a subalgebra of $A$.
\end{prop}
\begin{proof}
We first make some introductory remarks. Let $x = \sum_{j=1}^n r_j
\ot s_j$ and $y = \sum_{i=1}^m p_i \ot q_i$ be elements of $A \od
A$. Define
$$M(x,y)=\sum_{j,i=1}^{n,m} r_j p_i \ot q_i s_j.$$
Then we have
\begin{align*}
\|M(x,y)\|_h^2 & \leq \Bigl\|\sum_{i,j} r_j p_i p_i^* r_j^* \Bigr\| \, \Bigl\|
\sum_{i,j} s_j^* q_i^* q_i s_j \Bigr\| \\
& \leq \Bigl\|\sum_{i} p_i p_i^* \Bigr \| \, \Bigl\|\sum_{j} r_j r_j^* \Bigr\| \,
\Bigl\| \sum_{i} q_i^* q_i \Bigr\| \, \Bigl\| \sum_{j} s_j^* s_j \Bigr\|.
\end{align*}
From this we can conclude that $\| M(x,y) \|_h \leq \|x\|_h \,
\|y\|_h$, and we extend $M$ to a contraction from $(A \oh A)
\times (A \oh A)$ to $A \oh A$ by continuity. If now $y \in A \omh A$
and $x = \sum_{j =1}^n r_j \ot s_j \in A \od A$, we define,
consistently with the previous notation if $y \in A \oh A$
$$M(x,y)= \sum_{j = 1}^n (r_j \ot 1) y (1 \ot s_j) \in A \oh A.$$
Then we have, with $(e_\lambda)$ an approximate identity for $A$
$$\|M(x,y)\|_h = \lim_{\lambda} \Bigl\| \sum_{j=1}^n (r_j \ot 1)y(1 \ot
e_\lambda s_j) \Bigr\|_h = \lim_\lambda \|M(x,y(1 \ot e_\lambda))\|_h
\leq \|x\|_h \, \|y\|_h.$$
So again we can extend $M$ to a contraction from $(A \oh A) \times (A \omh
A)$ to $A \oh A$. If finally $x, y \in A \omh A$ we define $z \in
A \omh A$ by
\begin{equation*}
z(1 \ot c) = M(x(1 \ot c), y) \quad\quad
(c \ot 1) z = M((c \ot 1)x,y)
\end{equation*}
for all $c \in A$.
Suppose now $a,b \in A_0$. Take $x,y \in A \omh A$ satisfying
$T_1(x) = a \ot 1$ and $T_1(y)= b \ot 1$. Put $z = M(x,y)$. Then
we have $z \in A \omh A$ and we compute $T_1(z)$. Let $c \in A$.
Then
$$T_1(z)(1 \ot c)=T_1(z (1 \ot c)) = T_1(M(x(1 \ot c),y)).$$
If $p,q \in A$ we have
\begin{multline*}
T_1(M(p \ot q, y)) = T_1((p \ot 1) y (1 \ot q)) \\
= \phi(p) T_1(y) (1 \ot q) = \phi(p) (1 \ot q) (b \ot 1) = T_1(p
\ot q) (b \ot 1).
\end{multline*}
By continuity we may conclude from this and the previous formula
that
$$T_1(z) (1 \ot c) = T_1(x(1 \ot c))(b \ot 1) = (a \ot c)(b \ot
1)= ab \ot c.$$
So we get $T_1(z)= ab \ot 1$ and hence $a b \in A_0$.
\end{proof}

Next, we need to work with the map
$$\psi : p \otimes q \mapsto \phi_{13}(p) \phi_{23}(q)$$
from $A \otimes A$ to $M(A \otimes A \otimes A)$ where we use the common leg
numbering notation : $\phi_{13}(p)$ stands for $\phi(p)$ sitting
in the first and third factor, while $\phi_{23}(q)$ is simply $1 \otimes
\phi(p)$.  We first extend this map to $A \otimes_h A$ :

\begin{prop} \label{33}
The map $\psi$ is continuous for the Haagerup norm and
extends to a contraction from $A \otimes_h A$ to $M(A \otimes A \otimes A)$.
\end{prop}
\begin{proof}
This is analogous to the proof of proposition~\ref{13}. Just
remark that by lemma~\ref{11b}
\begin{align*}
\Bigl\|\psi \Bigl( \sum_{i=1}^n p_i \ot q_i \Bigr) \Bigr\|^2 &= \Bigl\| \sum_{i=1}^n
\phi_{13}(p_i) \phi_{23}(q_i) \Bigr\|^2 \\
& \leq \Bigl\| \sum_{i=1}^n \phi_{13}(p_i p_i^*) \Bigr\| \,
\Bigl\| \sum_{i=1}^n \phi_{23}(q_i^* q_i) \Bigr\| \\
& \leq \Bigl\| \sum_{i=1}^n p_i p_i^* \Bigr\| \,
\Bigl\| \sum_{i=1}^n q_i^* q_i \Bigr\|.
\end{align*}
So we get $\|\psi (x) \| \leq \|x\|_h$ for all $x
\in A \od A$ and that proves the result.
\end{proof}
We can say more about this map, but before we give the proof of
the next proposition we need the following lemma.
\begin{lemma} \label{strict}
Let $x \in M(A \ot A)$ such that $(\io \ot \om)(x) \in A$ and $(\io \ot \om)(x x^*) \in A$
for all $\om \in A^*$. Let $A$ act on a Hilbert space $H$, $\xi \in H$ and
let $(\eta_i)_{i \in I}$ be an orthonormal basis for $H$. Then the
increasing net in $A^+$
\begin{equation} \label{stijgnet}
\sum_{i \in I} (\io \ot \om_{\eta_i,\xi})(x) (\io \ot
\om_{\eta_i,\xi})(x)^*
\end{equation}
is norm convergent in $A$ with limit $(\io \ot
\om_{\xi,\xi})(xx^*)$.
\end{lemma}
Remark that we use the notation $\om_{\xi,\eta}$ for the vector
functional given by $\om_{\xi,\eta}(a) = \langle a \xi,\eta
\rangle$.
\begin{proof}
Choose $\om \in A^*_+$. Let $\pi$ be the universal representation of
$A$ on the Hilbert space $K$. Take $\rho\in K$ such that $\om(x)=
\om_{\rho,\rho}(\pi(x))$ for all $x \in A$. Let $(\rho_j)_{j \in
J}$ be an orthonormal basis for $K$. Let $I_0$ be a finite subset
of $I$. Then we have
\begin{align*}
\om \Bigl( \sum_{i \in I_0} (\io \ot \om_{\eta_i,\xi})(x) (\io \ot
\om_{\eta_i,\xi})(x)^* \Bigr) &= \sum_{i \in I_0} \| (\io \ot
\om_{\xi,\eta_i}) (\pi \ot \io)(x^*) \rho \|^2 \\
&= \sum_{i \in I_0, j \in J} | \langle (\pi \ot \io)(x^*) (\rho
\ot \xi), \rho_j \ot \eta_i \rangle |^2.
\end{align*}
So the increasing net of positive numbers
$$\Bigl(\om \Bigl( \sum_{i \in I_0} (\io \ot \om_{\eta_i,\xi})(x) (\io \ot
\om_{\eta_i,\xi})(x)^* \Bigr) \Bigr)_{I_0 \subseteq I}$$
converges to
$$\|(\pi \ot \io)(x^*) (\rho \ot \xi) \|^2 = \om ((\io \ot
\om_{\xi,\xi})(xx^*)).$$
Now we can apply Dini's theorem to the compact Hausdorff space $\{
\om \in A^*_+ \mid \|\om\| \leq 1 \}$, and the increasing
pointwise converging net of positive continuous functions on it given
by formula~\ref{stijgnet}. So we get uniform convergence of this net and
hence the lemma follows.
\end{proof}

\begin{prop} \label{34}
For all $\omega \in A^\ast$ and $x \in A \otimes_h A$
we have that $(\iota \otimes \iota \otimes \omega) \psi(x) \in A \otimes_h A$
and
$$\| (\iota \otimes \iota \otimes \omega)\psi(x)\|_h \leq \| \omega \| \, \| x
\|_h.$$
\end{prop}
\begin{proof}
By continuity it is enough to verify this for $x \in A \od A$. So
let $x = \sum_{i=1}^n p_i \ot q_i$ and choose $\om \in A^*$.
Let $A$ act on its universal Hilbert space $H$ and take vectors $\xi$
and $\eta$ in $H$ such that $\om = \om_{\xi,\eta}$ and
$\|\om\|=\|\xi\| \, \|\eta\|$. Let $(\gamma_j)_{j \in J}$ be an
orthonormal base for $H$. By the previous lemma we get that for
every $i=1,\ldots,n$ the net
$$\sum_{j \in J} ( \iota \otimes \omega_{\gamma_j,\eta}) \phi (p_i) (\iota \otimes
\omega_{\gamma_j,\eta}) \phi(p_i)^\ast$$
is norm convergent in $A$ with limit $(\io \ot
\om_{\eta,\eta})\phi(p_i p_i^*)$. Using the remark after
definition~\ref{11} we may conclude from this and an analogous
statement for the second leg that
\begin{equation} \label{vergelijking}
\sum_{i,j} (\iota \otimes
\omega_{\gamma_{j,\eta}}) \phi(p_i) \otimes (\iota \otimes
\omega_{\xi,\gamma_j})\phi (q_i) \in A \oh A
\end{equation}
with norm convergence. The norm of this element is majorized by
$$\Bigl\| \sum_{i=1}^n p_i p_i^* \Bigr\|^{1/2} \, \Bigl\| \sum_{i =1}^n q_i^* q_i
\Bigr\|^{1/2} \, \|\xi\| \, \| \eta \|.$$
Because the net in~\ref{vergelijking} is also norm convergent in
$A \ot A$ with limit
$$\sum_{i=1}^n (\io \ot \io \ot \om_{\xi,\eta}) (\phi_{13}(p_i)
\phi_{23}(q_i)) = (\io \ot \io \ot \om) \psi(x)$$
the result follows.
\end{proof}

Now we extend $\psi$ to $A \otimes_{mh}A$ by strict continuity:

\begin{prop} \label{35}
The map $\psi$ can be extended uniquely to a contraction from $A
\omh A$ to $M(A \ot A \ot A)$, still denoted by $\psi$, such that
for all $x \in A \omh A$ and $a \in A$
\begin{equation*}
\psi(x) \phi_{23}(a) = \psi(x(1 \ot a)) \quad\text{and}\quad
\phi_{13}(a) \psi(x) = \psi((a \ot 1)x).
\end{equation*}
This extension is strictly continuous on bounded sets. For every
$\om \in A^*$ the map
$$ x \mapsto (\io \ot \io \ot \om) \psi(x)$$
maps $A \omh A$ into $A \omh A$, is bounded with norm majorized by
$\|\om\|$ and is strictly continuous on bounded sets.
\end{prop}
\begin{proof}
The extension of $\psi$ to $A \omh A$ is completely analogous to
the extension of the maps $T_1$ and $T_2$. If $x \in A \omh A$ we can obtain $\psi(x)$
as a right multiplier by putting $(a \ot b \ot c) \phi_{13}(d)
\psi(x)$ equal to $(a \ot b \ot c) \psi((d \ot 1)x)$. One shows
just as in the proof of proposition~\ref{25} that we obtain a
bounded right multiplier of $A \ot A \ot A$ this way. In a similar
way we define a bounded left multiplier. This gives us the
extension of $\psi$.

Now write for every $\om \in A^*$, $x \in A \oh A$ and $u \in A
\ot A$
$$R(\om,x,u) = (\io \ot \io \ot \om)(u_{13} \psi(x)).$$
At the moment we only have $R(\om,x,u) \in M(A \ot A)$. We claim
that $R(\om,x,u) \in A \oh A$ and
\begin{equation} \label{vglR}
\|R(\om,x,u)\|_h \leq \|\om\| \, \|x\|_h \, \|u\|.
\end{equation}
When $x = \sum_{i=1}^n p_i \ot q_i$, when $A$ acts on a Hilbert
space $H$ with $\om = \om_{\xi,\eta}$, $\|\om\|=\|\xi\|\,\|\eta\|$
and when $(\gamma_j)_{j \in J}$ is an orthonormal basis for $H$,
then we have as in the proof of the previous proposition
$$R(\om,x,u)= \sum_{i,j} (\io \ot \om_{\gamma_j,\eta})(u
\phi(p_i)) \ot (\io \ot \om_{\xi,\gamma_j})(\phi(q_i))$$
with norm convergence in $A \oh A$. Further we get
$$\|R(\om,x,u)\|_h  \leq \|\eta\| \, \|u\| \, \Bigl\| \sum_{i=1}^n
p_i p_i^* \Bigr\|^{1/2} \, \|\xi\| \, \Bigl\| \sum_{i=1}^n q_i^* q_i \Bigr\|^{1/2}.$$
This proves our claim that $R(\om,x,u) \in A \oh A$ and also
proves inequality~\ref{vglR}.

Now let $\om \in A^*$, $x \in A \omh A$, $u \in A \ot A$ and $a
\in A$. Then
$$(\io \ot \io \ot \om)((u \phi(a))_{13} \psi(x)) = (\io \ot \io
\ot \om)(u_{13} \psi((a \ot 1) x)) \in A \oh A$$
and
\begin{align*}
\|(\io \ot \io \ot \om)((u \phi(a))_{13} \psi(x))\|_h &=
\|R(\om,(a \ot 1) x, u)\|_h \\
&= \lim_\lambda \|R(\om,(a \ot 1) x (1 \ot e_\lambda),u)\|_h \\
&= \lim_\lambda \|R(\om,x(1 \ot e_\lambda),u\phi(a)) \|_h \\
&\leq \|\om\| \, \|x\|_h \, \|u \phi(a)\|.
\end{align*}
Because $\phi$ is non-degenerate we may conclude that for $\om \in
A^*$, $x \in A \omh A$ and $u \in A \ot A$ we have $(\io \ot \io
\ot \om)(u_{13} \psi(x)) \in A \oh A$. We still denote this
element with $R(\om,x,u)$ and we have
$$\|R(\om,x,u)\|_h \leq \|\om\| \, \|x\|_h \, \|u\|.$$
From the formulas above we can also conclude that for every $\om
\in A^*$ and $u \in A \ot A$ the map
$$ x \mapsto R(\om,x,u)$$
is strict-norm continuous on bounded sets from $A \omh A$ to $A
\oh A$.

Finally let $\om \in A^*$ and $x \in A \omh A$. Write $\om = \mu b$
with $\mu \in A^*$ and $b \in A$. Then we get for every $a \in A$
$$(a \ot 1) (\io \ot \io \ot \om)\psi(x) = (\io \ot \io \ot
\mu)((a \ot b)_{13} \psi(x)) \in A \oh A$$
and
\begin{align*}
\|(a \ot 1) (\io \ot \io \ot \om)\psi(x)\|_h &= \|R(\mu,x, a \ot
b)\|_h \\
&= \lim_\lambda \| R(\mu,x,a \ot be_\lambda) \|_h \\
&= \lim_\lambda \| R(\om,x,a \ot e_\lambda)\|_h \\
& \leq \|\om\| \, \|x\|_h \, \|a\|.
\end{align*}
Analogously one can prove that $(\io \ot \io \ot \om)\psi(x)(1
\ot a) \in A \oh A$ with norm majorized by $\|\om\| \, \|x\|_h \,
\|a\|$. So we get indeed $(\io \ot \io \ot \om)\psi(x) \in A
\omh A$ with norm majorized by $\|\om\| \, \|x\|_h$. With the
notation as above we can write
$$(a \ot 1) (\io \ot \io \ot \om)\psi(x) = R(\mu,x,a \ot b)$$
and this depends strict-norm continuous on bounded $x$ by one of the
remarks above. One proves an analogous statement on the other side
and gets the strict continuity on bounded sets of $x \mapsto (\io \ot \io \ot
\om)\psi(x)$ from $A \omh A$ to $A \omh A$.
\end{proof}

These results allow us to prove the following basic formula.

\begin{prop} \label{36}
If $a \in A_0$ and $x \in A \otimes_{mh} A$ such that
$a \otimes 1 = T_1 x$, then $x \otimes 1 = \psi (x)$.
\end{prop}
\begin{proof}
We will prove this formula by applying $\omega$ and by showing
that $T_1((\iota \otimes \iota \otimes \omega)\psi(x)) = T_1(x) \omega (1)$.
Then, the injectivity of $T_1$ will give us the result.

Now observe that
\begin{align*}
T_1 (( \iota \otimes \iota \otimes \omega)\psi(p \otimes q)) & = T_1 ((\iota
\otimes \iota \otimes \omega) (\phi_{13} (p) \phi_{23}(q)) )\\
& = (\iota \otimes \iota \otimes \omega)((\phi \otimes \iota)\phi(p) \phi_{23}
(q)).
\end{align*}
This can easily be shown by using the concrete formulas for
$\psi(p \ot q)$ obtained in the proof of proposition~\ref{34}.
Then, using coassociativity, we get
\begin{align*}
T_1((\iota \otimes \iota \otimes \omega) \psi(p \otimes q)) & = (\iota \otimes
\iota \otimes \omega) ((\iota \otimes \phi)(\phi(p)(1 \otimes q))) \\
& = (\iota \otimes \iota \otimes \omega)((\iota \otimes \phi)T_1(p \otimes q)).
\end{align*}
Then, using the strict continuity on bounded sets of all the maps involved, we get
$$T_1((\iota \otimes \iota \otimes \omega)\psi(x)) = (\iota \otimes \iota
\otimes \omega)((\iota \otimes \phi)T_1(x)).$$
As $T_1x = a \otimes 1$, the righthand side is
\begin{align*}
(\iota \otimes \iota \otimes \omega)((\iota \otimes \phi)(a  \otimes 1)) & = (a
\otimes 1) \omega (1) \\
& = T_1(x) \omega (1).
\end{align*}
Then the proof is complete.
\end{proof}

Before we continue, let us make a little excursion to Hopf algebras and see what
the above formula really means.  We will use the Sweedler notation.  So, let
$(H,\Delta, \varepsilon,S)$ be a Hopf algebra.  For any $a \in H$ we have
\begin{align*}
a \otimes 1 & = \sum a_{(1)} \otimes \varepsilon(a_{(2)})1 = \sum a_{(1)} \otimes a_{(2)}
S(a_{(3)}) \\
& = \sum \Delta (a_{(1)})(1 \otimes S(a_{(2)})).
\end{align*}
So, in this case $a \otimes 1 = T_1 (\sum a_{(1)} \otimes S(a_{(2)}))$.  Now
\begin{align*}
\sum \Delta_{13} (a_{(1)}) \Delta_{23} (S(a_{(2)})) & = \sum (a_{(1)} \otimes 1
\otimes a_{(2)})(1 \otimes S(a_{(4)}) \otimes S(a_{(3)})) \\
& = \sum a_{(1)} \otimes S(a_{(4)}) \otimes a_{(2)} S(a_{(3)}) \\
& = \sum a_{(1)} \otimes S(a_{(3)}) \otimes \varepsilon (a_{(2)})1 \\
& = \sum a_{(1)} \otimes S(a_{(2)}) \otimes 1
\end{align*}
and this is precisely the formula that we have, in the C$^*$-context, in
proposition~\ref{36}.

This is very important for understanding what follows.  Simply observe that
$$\sum a_{(1)} S(a_{(2)}) = \varepsilon(a)1$$
and
$$ \sum (1 \otimes a_{(1)}) \Delta (S(a_{(2)})) = \sum
S(a_{(3)}) \otimes a_{(1)} S(a_{(2)})
 = S(a) \otimes 1.
$$
These formulas will motivate the definition of $\varepsilon$ and $S$ below.

First we need a little lemma.
\begin{lemma} \label{37}
If $x \in M(A)$ and $\phi(x) = x \otimes 1$, then $x \in
\C 1$.
\end{lemma}
\begin{proof}
Take $a,b \in A$.  Then
$$(a \otimes 1) \phi (xb) = (a \otimes 1)\phi(x)\phi(b) = (ax \otimes
1)\phi(b).$$

We now use the injectivity of $T_2$ (on $A \odot A$) to find
$$a \otimes xb = ax \otimes b$$
for all $a,b \in A$.  This is only possible when $x \in \C 1$.
\end{proof}
This lemma, together with the formulas in the Hopf algebra context above will
lead us to the definition of the counit and the antipode.
\begin{prop} \label{38}
There is a homomorphism $\varepsilon : A_0 \rightarrow
\C$, defined by $\varepsilon(a) 1 = m(x)$ when $a \in A_0$ and $x \in A \otimes_{mh}
A$ satisfy $a \otimes 1 = T_1 x$.  Here $m$ denotes multiplication as a map
from $A \otimes_{mh} A$ to $M(A)$ (see proposition~\ref{19}).
\end{prop}
\begin{proof}
Take $a \in A_0$ and $x \in A \otimes_{mh} A$ such that $a
\otimes 1 = T_1 (x)$.  Take $\omega \in A^\ast$ and consider the formula
$\omega(1) x = (\iota \otimes \iota \otimes \omega)\psi(x)$.  Apply
multiplication on both sides.  For any $p,q \in A$ we have
\begin{align*}
m(\iota \otimes \iota \otimes \omega)(\phi_{13} (p) \phi_{23}(q)) & = (\iota
\otimes \omega) (\phi (p) \phi(q)) = (\iota \otimes \omega)\phi(pq) \\
& = (\iota \otimes \omega) \phi (m(p \otimes q)).
\end{align*}
Again this can be proved by writing out the formula for $(\iota \otimes \iota \otimes \omega)(\phi_{13} (p)
\phi_{23}(q))$ that we obtained in the proof of proposition~\ref{34}.
We get by continuity that
$$\omega(1) m(x) = m(\io \ot \io \ot \om) \psi(x) = (\iota \otimes \omega) \phi (m(x)).$$
As this is true for all $\omega$ we must have
$$m(x) \otimes 1 = \phi (m(x)).$$
By the lemma $m(x) \in \C 1$.  So we can define $\varepsilon : A_0 \rightarrow \C$
by $\varepsilon (a)1 = m(x)$ when $a \otimes 1 = T_1x$ and $x \in A \otimes_{mh}
A$.

Now we show that $\eps$ is a homomorphism.  To do this, take $a,b \in A_0$ and
$x,y \in A \otimes_{mh} A$ such that $T_1 x = a \otimes 1$ and $T_1 y = b
\otimes 1$. Use the notation of the proof of
proposition~\ref{32} and define $z= M(x,y)$. We have seen that $z
\in A \omh A$ and $T_1z = ab \ot 1$. Let $c \in A$. Then we have
$$\eps(ab)c = m(z(1 \ot c)) = m(M(x(1 \ot c),y)).$$
But for $p,q \in A$ we have
$$m(M(p \ot q, y)) = m((p \ot 1) y (1 \ot q)) = p m(y) q = pq
\eps(b) = m(p \ot q) \eps(b).$$
So we can conclude that
$$\eps(ab)c = m(x(1 \ot c)) \eps(b) = \eps(a) \eps(b) c.$$
This gives the required result.
\end{proof}

The antipode is obtained in a similar way.  Instead of $m$, we now need to use
the map $p \otimes q \rightarrow (1 \otimes p)\Delta(q)$.  This can be seen
from the motivation explained  before in the Hopf algebra context.  We will
call this map $\widetilde T$.  It is of the same kind as the maps $T_1$ and $T_2$
(with $p$ and $q$ in the correct order).  Therefore one can use the same
methods to extend $\widetilde T$ to a continuous map from $A \otimes_{mh} A$ to
$M(A \otimes A)$, continuous both for the norm and for the strict topology on
bounded sets.  It is also possible to view $\widetilde T$ as
$$\widetilde T(p \otimes q) = (T_1(q^\ast \otimes p^\ast)^\ast.$$
Now remark that $p \otimes q \mapsto q^\ast \otimes p^\ast$ extends to an
isometry from $A \otimes_h A$ to itself and further to an isometry from $A
\otimes_{mh} A$ to itself that is strictly continuous (see proposition~\ref{112}).

This takes us to the definition of the antipode.
\begin{prop} \label{39}
There is a linear map $S : A_0 \rightarrow A$ defined
by $S(a) \otimes 1 = \widetilde T(x)$ when $a \otimes 1 = T_1(x)$ and $a \in A_0$ and
$x \in A \otimes_{mh} A$.  This map is an anti-homomorphism.  We have $S(a)^\ast
\in A_0$ when $a \in A_0$ and $S(S(a)^\ast)^\ast = a$.
\end{prop}
\begin{proof}
Let $a \in A_0$ and $x \in A \otimes_{mh} A$ satisfy $a \otimes 1
= T_1 x$.  Apply the slice map $\iota \otimes \iota \otimes \omega$ and then
$\widetilde T$ to the formula $x \otimes 1 = \psi (x)$ as obtained in proposition~\ref{36}.
This yields
$$\omega(1) \widetilde T(x) = \widetilde T ( \iota \otimes \iota \otimes \omega) \psi
(x).$$
Now, when $p,q \in A$ we have with the same technique as in the proof of the previous proposition
\begin{align*}
\widetilde T(\iota \otimes \iota \otimes \omega) \psi (p \otimes q) & = \widetilde T
(\iota \otimes \iota \otimes \omega) (\phi_{13}(p) \phi_{23} (q)) \\
& = (\iota \otimes \iota \otimes \omega)(\phi_{23} (p) (\phi \otimes \iota)
\phi(q)) \\
& = (\iota \otimes \iota \otimes \omega)(\iota \otimes \phi)((1 \otimes p)
\phi(q)) \\
& = (\iota \otimes \iota \otimes \omega)(\iota \otimes \phi) \widetilde T(p \otimes
q).
\end{align*}
By continuity this is still true when we replace $p \otimes q$ by $x$.  We
obtain
$$\omega(1) \widetilde T(x) = \widetilde{T}(\io \ot \io \ot \om) \psi(x)= (\iota \otimes \iota \otimes \omega)(\iota \otimes
\phi) \widetilde T(x).$$
This gives  $\widetilde T(x) \otimes 1 = (\iota \otimes \phi) \widetilde T(x)$.  Now
apply the slice map $\omega \otimes \iota \otimes \iota$ to get
$$(\omega \otimes \iota) \widetilde T(x) \otimes 1 = \phi ((\omega \otimes  \iota)
\widetilde T(x)).$$
By lemma~\ref{37} we get $(\omega \otimes \iota) \widetilde T(x) \in \C 1$.  As this is
true for all $\omega$, we must have $\widetilde{T}(x) \in M(A) \otimes 1$.  As in the
remark following~\ref{31} we see that actually $\widetilde T(x) \in A \otimes 1$.  Then
we can define $S(a) \in A$ by $\widetilde T(x) = S(a) \otimes 1$ when $a \in A_0$ and $x
\in A \otimes_{mh} A$ satisfy $a = T_1 x$.

Next we show that $S$ is an anti-homomorphism.  As before take $a,b \in A_0$ and
$x,y \in A \otimes_{mh} A$ such that $a \otimes 1 = T_1x$ and $b \otimes 1  =
T_1 y$. Use again the notation of the proof of
proposition~\ref{32}. Put $z = M(x,y)$. Then $z \in A \omh A$ and
$T_1 z=ab \ot 1$. Let $c \in A$. Then we get
$$(S(ab) \ot 1) \phi(c) = \widetilde{T} (z(1 \ot c)) = \widetilde{T}(M(x(1
\ot c),y)).$$
Now for $p,q \in A$ one has
$$\widetilde{T}(M(p \ot q,y)) = \widetilde{T}((p \ot 1)y(1 \ot q)) = (1 \ot p) \widetilde{T}(y) \phi(q) =
(S(b) \ot 1) \widetilde{T}(p \ot q).$$
This gives
$$(S(ab) \ot 1)\phi(c) = (S(b) \ot 1) \widetilde{T}(x(1 \ot c)) =
(S(b)S(a) \ot 1) \phi(c)$$
and so $S(ab)=S(b) S(a)$. Finally
$$
S(a)^\ast \otimes 1  = (S(a) \otimes 1)^\ast
 = T_1 (\Gamma(x))
$$
where $\Gamma$ is the map from proposition~\ref{112}.  As $\Gamma(x) \in A \otimes_{mh} A$,
we have $S(a)^\ast \in A_0$.  Also
$$
S(S(a)^\ast) \otimes 1 = \widetilde T (\Gamma(x))
= T_1 (x)^\ast = a^\ast \otimes 1
$$
and this proves the result.
\end{proof}
In the case of a multiplier Hopf algebra we have essentially the
same formulas. The counit in a multiplier Hopf algebra has been
defined in just the same way as we do it here, see section~3 of
\cite{VD2}. The way to define the antipode for multiplier Hopf
algebras was not usable however.  The formulas in
lemma~5.4 and lemma~5.5 of~\cite{VD2} are nevertheless very similar to the
formulas we use to define the antipode.

We would now like to prove some of the well-known formulas involving the counit
and the antipode.  We will refer to the discussion above in the Hopf algebra
context.  Recall that
$$
a \otimes 1 = \sum \Delta (a_{(1)})(1 \otimes S(a_{(2)}))
= T_1 ((\iota \otimes S) \Delta (a)).
$$
So we expect that the element $x \in A \otimes_{mh} A$ satisfying $T_1 (x) = a
\otimes 1$ should be, in some sense, $(\iota \otimes S)\Delta(a)$.  This is the
content of the following proposition.

\begin{prop} \label{310}
Let $a \in A_0, x \in A \otimes_{mh} A$ and $a
\otimes 1 = T_1 x$.  For all $\omega \in A^\ast$ we have that $(\omega \otimes
\iota) \phi (a) \in A_0$ and $S((\omega \otimes \iota)\phi(a)) = (\omega
\otimes \iota)(x)$.
\end{prop}
\begin{proof}
Before we can really start the proof we have to make some
introductory remarks. Let $p,q \in A$ and $\om \in A^*$. Then we
have
\begin{align*}
(\omega \otimes \iota \otimes \iota)\psi (p \otimes q) & = (\omega \otimes
\iota \otimes \iota)(\phi_{13}(p) \phi_{23} (q)) \\
& = (1 \otimes (\omega \otimes \iota) \phi(p)) \phi(q) \\
& = \widetilde T ((\omega \otimes \iota)\phi(p) \otimes q) \\
& = \widetilde T (((\omega \otimes \iota)\phi \otimes \iota)(p \otimes q)).
\end{align*}
Now we would like to have the same formula for arbitrary $x \in A
\omh A$ instead of $p \ot q$. So we have to extend the map $(\om
\ot \io)\phi \ot \io$. Because $\om \ot \io$ is completely bounded
as a map from $M(A \ot A)$ to $M(A)$ with norm $\|\om\|$ we get
that $(\om \ot \io)\phi \ot \io$ can be uniquely extended to a
continuous map from $A \oh A$ to itself with norm majorized by
$\|\om \|$. See the discussion in the beginning of section~4.
For every $\om \in A^*$ and $u \in A \ot A$ we can, by
the same argument, extend
$$p \ot q \mapsto (\om \ot \io)(u \phi(p)) \ot q$$
to a continuous map from $A \oh A$ to itself with norm majorized
by $\|\om\|\,\|u\|$. We write
$$Q(\om,u,x) = ((\om \ot \io)(u \phi( \cdot)) \ot \io)(x)$$
for this extension applied to $x \in A \oh A$. If now $x
\in A \omh A$, $u \in A \ot A$, $\om \in A^*$ and $a \in A$ we get
\begin{align*}
\|Q(\om,u,(a \ot 1)x)\|_h &= \lim_\lambda \| Q(\om,u,(a \ot 1) x
(1 \ot e_\lambda)) \|_h \\
&= \lim_\lambda \| Q(\om,u\phi(a), x(1 \ot e_\lambda)) \|_h \\
&\leq \|\om\| \, \|u \phi(a) \| \, \|x\|_h
\end{align*}
So by the non-degenerateness of $\phi$ it is possible to extend
$Q$ uniquely to a map
$$A^* \times (A \ot A) \times (A \omh A) \rightarrow A \oh A$$
such that
$$Q(\om,u\phi(a),x) = Q(\om, u , (a \ot 1) x)$$
for $\om \in A^*$, $u \in A \ot A$, $a \in A$ and $x \in A \omh
A$. It is clear from the formulas above that for every $\om \in A^*$ and $u \in
A \ot A$ the map $x \mapsto Q(\om,u,x)$ from $A \omh A$ to $A \oh
A$ is strict-norm continuous on bounded sets.

Let now $\om \in A^*$. Write $\om = \mu b$
with $\mu \in A^*$ and $b \in A$. Then define
$$\rho_2(a,x)=Q(\mu,b \ot a,x) \in A \oh A$$
for every $a \in A$ and $x \in A \omh A$. Then for every $a \in A$
the map $x \mapsto \rho_2(a,x)$ is strict-norm continuous on
bounded sets and
$$\rho_2(a,p \ot q) = (a \ot 1) ((\om \ot \io)\phi \ot \io)(p \ot
q).$$
Analogously the map
$$p \ot q \mapsto ((\om \ot \io) \phi \ot \io)(p \ot q) (1 \ot
a)$$
can be extended uniquely to a map from $A \omh A$ to $A \oh A$
which is strict-norm continuous on bounded sets. So we have proved
that $(\om \ot \io)\phi \ot \io$ can be extended to a bounded map
from $A \omh A$ to itself which is strictly continuous on bounded
sets. We denote this (unique) extension by the same symbols.

Now we start the real proof. Let $a \in A_0$, $x \in A \omh A$
and $T_1x = a \ot 1$. Then we have $x \ot 1 = \psi(x)$ so that for
all $\om \in A^*$ we get $(\om \ot \io)(x) \ot 1 = (\om \ot \io
\ot \io)\psi(x)$. Now the computation in the beginning of this
proof, together with the strict continuity of the maps involved
gives
$$(\omega \otimes \iota)(x) \otimes 1 =(\om \ot \io \ot \io) \psi(x)= \widetilde T(((\omega \otimes \iota)\phi
\otimes \iota)(x)).$$
It follows from this that $(\omega \otimes \iota)(x) = S(b)$ for $b \in A_0$
satisfying
$$b \otimes 1 = T_1(((\omega \otimes \iota) \phi \otimes \iota)(x)).$$
Now
\begin{align*}
T_1(( (\omega \otimes \iota)\phi \otimes \iota)(p \otimes q)) & = (\omega \otimes
\iota \otimes \iota) ((\iota \otimes \phi)\phi(p)(1 \otimes 1 \otimes q)) \\
& = ((\omega \otimes \iota) \phi \otimes \iota)(T_1(p \otimes q)).
\end{align*}
So if we replace $p \otimes q$ by $x$ we get
$$
b \otimes 1  = ((\omega \otimes \iota)\phi \otimes \iota) (T_1(x))
= ((\omega \otimes \iota) \phi \otimes \iota) (a \otimes 1)
$$
and so
$b = (\omega \otimes \iota)\phi(a)$.
\end{proof}

So formally $x = (\iota \otimes S)\phi(a)$ when $a \otimes 1 = T_1(x)$.  Now,
by definition $\varepsilon(a)1 = m(x)$ and we recover the formula $m(\iota \otimes
S) \phi (a) = \varepsilon(a)1$.  Actually we can now give a meaning to this
formula because we take elements $a$ such that $(\iota \otimes S)\Delta(a)$
belongs to $A \otimes_{mh}A$, where we can apply multiplication (recall the
discussion in the introduction).

From the proof of the previous proposition we also find $\varepsilon((\omega \otimes \iota)\phi(a))1 = m(((\omega \otimes
\iota)\phi \otimes \iota)(x))$, when $a \otimes 1 = T_1 x$ and $\omega \in
A^\ast$.  Now
$$
m(((\omega \otimes \iota)\phi \otimes \iota)(p \otimes q))  = (\omega \otimes
\iota)(\phi(p)(1 \otimes q)) = (\omega \otimes \iota)T_1(p \otimes q).
$$
If we use this formula for $x$ in the place of $p \otimes q$, we arrive at
$$\varepsilon((\omega \otimes \iota)\phi(a))1 = \omega (a) 1$$
which precisely stands for $(\iota \otimes \varepsilon)\phi(a) = a$ when $a \in
A_0$.

Also the following formula is very familiar:
\begin{prop}
Let $a \in A_0$. Then $S(a)^* \in A_0$ and $\eps(S(a)^*) =
\overline{\eps(a)}$.
\end{prop}
\begin{proof}
Take $x \in A \omh A$ such that $T_1(x)=a \ot 1$. It follows from
the proof of proposition~\ref{39} that $S(a)^* \ot 1 =
T_1(\Gamma(x))$. Hence we have $\eps(S(a)^*)1 = m(\Gamma(x))$. Now
for $p,q \in A$ we get
$$m(\Gamma(p \ot q)) = m(q^* \ot p^*) = q^* p^* = (m(p \ot q))^*$$
so that $m(\Gamma(x)) = (m(x))^* = \overline{\eps(a)}1$. This
gives the result.
\end{proof}

There is one more formula in Hopf algebra theory that we can give a meaning
here, namely
$$\Delta(S(a)) = \sigma (S \otimes S)\Delta(a).$$

\begin{prop} \label{311}
If $a \in A_0$ and $x \in A \otimes_{mh} A$ such that
$a \otimes 1 = T_1 x$, then $(\io \ot \om)(x) \in A_0$ for all $\omega
\in A^\ast$ and $S((\io \ot \om)(x)) = (\omega \otimes \iota)\phi(S(a))$.
\end{prop}
The formula mentioned before the proposition is obtained by
writing formally that $x =(\io \ot S) \phi(a)$.
\begin{proof}
Start again from the basic formula $x \otimes 1 = \psi(x)$ and
now apply the slice map $\iota \otimes \omega \otimes \iota$.  For $p,q \in A$
we have
\begin{align*}
(\iota \otimes \omega \otimes \iota)(\phi_{13}(p)\phi_{23}(q)) & = \phi(p)(1
\otimes (\omega \otimes \iota)\phi(q)) \\
& = T_1 (\iota \otimes (\omega \otimes \iota)\phi)(p \otimes q).
\end{align*}
Exactly as in the proof of the previous proposition we can extend
$\io \ot (\om \ot \io)\phi$ to a map from $A \omh A$ to itself
which is strictly continuous on bounded sets. So applying the
formula above to $x$ in the place of $p \otimes q$ we get
$$
(\io \ot \om)(x) \otimes 1  = (\iota \otimes \omega \otimes
\iota)\psi(x)
= T_1( (\iota \otimes (\omega \otimes \iota)\phi)(x)).
$$
Hence $(\io \ot \om)(x) \in A_0$ and
$$S((\io \ot \om)(x)) \otimes 1 = \widetilde T((\iota \otimes (\omega \otimes
\iota) \phi)(x)).$$
Again, with $x$ replaced by $p \otimes q$ we get
\begin{align*}
\widetilde T((\iota \otimes (\omega \otimes \iota) \phi)(p \otimes q)) & =
(\omega \otimes \iota \otimes \iota)((1 \otimes 1 \otimes p)
(\iota \otimes \phi) \phi(q)) \\
& = ((\omega \otimes \iota)\phi \otimes \iota) \widetilde T(p \otimes q).
\end{align*}
So
\begin{align*}
S((\io \ot \om)(x)) \otimes 1 & = ((\omega \otimes \iota) \phi \otimes
\iota)(\widetilde T(x)) \\
& = (\omega \otimes \iota)\phi(S(a)) \otimes 1.
\end{align*}
This proves the result.
\end{proof}

Again we also get
$$\varepsilon((\io \ot \om)(x))1 = m((\iota \otimes (\omega \otimes
\iota)\phi)(x)).$$
With $x$ replaced by $p \otimes q$ we get
\begin{equation*}
m((\iota \otimes (\omega \otimes \iota)\phi)(p \otimes q))  = (\omega \otimes
\iota)((1 \otimes p)\phi(q))
 = (\omega \otimes \iota) \widetilde T(p \otimes q)
\end{equation*}
and so
$$\varepsilon((\io \ot \om)(x))1 = (\omega \otimes \iota)(S(a) \otimes 1).$$
Thus $\varepsilon((\io \ot \om)(x)) = \omega (S(a))$.  This essentially is the
formula
$$
(\varepsilon \otimes \iota)((\iota \otimes S)\Delta (a)) = S(a).$$

\begin{remark} \label{remark}
In this section we have worked with $T_1$ and $\widetilde{T}$ to
define $A_0$, $S$ and $\eps$. It is for reasons of symmetry clear
that we could also work with $A'_0$, consisting of those elements
$a \in A$ such that $1 \ot a$ belongs to the range of $T_2$ on $A
\omh A$. We can then extend the map $p \ot q \mapsto \phi(p)(q \ot
1)$ to a map $\widehat{T}$ from $A \omh A$ to $M(A \ot A)$ and define
$S'$ and $\eps'$ on $A'_0$ such that
$$1 \ot S'(a) = \widehat{T}(x) \quad \quad \eps'(a) 1 = m(x)$$
when $1 \ot a = T_2(x)$. In the next section and in section~5 we
will describe natural subspaces of $A_0 \cap A'_0$ on which $S,S'$
and $\eps,\eps'$ coincide.
\end{remark}

\section{The two-sided counit and antipode}
When $A$ is a C$^*$-algebra there is a natural way to norm the
linear space $M_{nm}(A)$ of $n$ by $m$ matrices over $A$. When
$A$ is a closed $^*$-subalgebra of $B(H)$ this norm is obtained by
looking at $M_{nm}(A) \subseteq B(H^m,H^n)$, where $H^m$ denotes
the $m$-fold direct sum of $H$. Now when $(p_i) \in M_{1n}(A)$
and $(q_i) \in M_{n1}(A)$ we have
$$\|(p_i)\| = \Bigl\| \sum_{i=1}^n p_i p_i^* \Bigr\|^{1/2}
\quad\text{and}\quad \|(q_i)\| = \Bigl\| \sum_{i=1}^n q_i^* q_i
\Bigr\|^{1/2}.$$
Compared to definition~\ref{11} the following definition should
not be too surprising:

\begin{defin}
Let $A$ be a C$^*$-algebra and $x \in A \od A \od A$. Then we
define
\begin{multline*}
\|x\|_h = \inf \Bigl\{ \|(p_i)\| \, \|(q_{ij}) \| \, \| (r_j) \| \mid
(p_i) \in M_{1n}(A), (q_{ij}) \in M_{nm}(A), (r_j) \in
M_{m1}(A) \\ \text{and} \quad x = \sum_{i,j=1}^{n,m} p_i \ot
q_{ij} \ot r_j \Bigr\}.
\end{multline*}
\end{defin}
So again we take the infimum over all such representations of $x$
in the algebraic tensor product and one can again prove that this
is a norm on $A \od A \od A$. The completion for this norm will be
denoted by $A \oh A \oh A$.

Also here we will regularly want to extend maps defined on the
algebraic tensor product to the triple Haagerup tensor product.
The following obvious lemma is crucial.
\begin{lemma} \label{41b}
Let $B$ be a C$^*$-algebra, $(p_i) \in M_{1n}(B)$, $(q_{ij}) \in
M_{nm}(B)$ and $(r_j) \in M_{m1}(B)$. Then
$$\Bigl\| \sum_{i,j=1}^{n,m} p_i q_{ij} r_j \Bigr\| \leq \|(p_i)\| \, \|(q_{ij}) \| \, \| (r_j)
\|.$$
\end{lemma}
\begin{proof}
This follows immediately when one takes $B \subseteq B(H)$ and
$M_{kl}(B) \subseteq B(H^k,H^l)$.
\end{proof}
We recall that a map $F$ from a C$^*$-algebra $A$ to a
C$^*$-algebra $B$ is called completely bounded when there exists a
positive number $M$ such that $\|(F(a_{ij}))\| \leq M \|(a_{ij})\|$
for all $(a_{ij}) \in M_{nm}(A)$ and all $n,m \in \N$. The
smallest such number $M$ is denoted by $\|F\|_{cb}$. When
$\|F\|_{cb} \leq 1$ we say $F$ is a complete contraction.
Completely bounded maps can be tensored freely on the Haagerup
tensor product. For instance, when $F,G,H : A \rightarrow B$ are
completely bounded, it is clear nearly by definition that $F \ot G
\ot H$ defined from $A \od A \od A$ to $B \od B \od B$ is bounded
for the Haagerup norm. We also mention that every $\om \in A^*$ is
completely bounded when considered as a map from $A$ to $\C$ and
$\|\om\|_{cb} = \|\om \|$. So it is clear how to define slice maps
$(\io \ot \io \ot \om)$ or $(\om \ot \io \ot \io)$ from $A \oh A
\oh A$ to $A \oh A$. Finally it is also possible to define norms
on the spaces of matrices over the Haagerup tensor product. If for
instance $x \in M_{nm}(A \od A)$ we define
\begin{multline*}
\|x \|_h = \inf \Bigl\{ \|(p_{ik})\| \, \|(q_{kj})\| \mid r \in \N,
(p_{ik}) \in M_{nr}(A), (q_{kj}) \in M_{rm}(A) \\ \text{and}\quad
x_{ij}= \sum_{k=1}^r p_{ik} \ot q_{kj} \; \text{for all} \; i,j
\Bigr\}.
\end{multline*}
So it also makes sense to speak about completely bounded maps from
$A \oh A$ to $A$. When $F$ is such a map, of course $F \ot \io$
will again be a completely bounded map from $A \oh A \oh A$ to $A
\oh A$. All this is clearly worked out in
the theory of \emph{operator spaces} and we refer to
e.g.~\cite{B-P}, \cite{B-P2} and \cite{E-R}
for this.

We will now treat the two-sided approach somewhat different from the one-sided
approach.  There we worked with the multiplier Haagerup tensor product $A
\otimes_{mh} A$ which could be characterized as the space of elements $x \in
M(A \otimes A)$ such that both $(a \otimes 1)x$ and $x(1 \otimes a)$ belong to
$A \otimes_h A$ for all $a \in A$.  We could also proceed here in a similar
way.  However, we will stick to the triple Haagerup tensor product.

First we need the map
$$T : p \otimes q \otimes r \mapsto (p \otimes 1 \otimes 1) \phi^{(2)}(q)(1
\otimes 1 \otimes r)$$
where $\phi^{(2)}(q) = (\phi \otimes \iota)\phi(q) = (\iota \otimes
\phi)\phi(q)$.  Remark that
\begin{align*}
T(p \otimes q \otimes r) & = (\iota \otimes \phi)((p \otimes 1)\phi(q))(1
\otimes 1 \otimes r) \\
T(p \otimes q \otimes r) & = (p \otimes 1 \otimes 1)(\phi \otimes
\iota)(\phi(q)(1 \otimes r))
\end{align*}
so that formally $T = (\iota \otimes T_1)(T_2 \otimes \iota) = (T_2 \otimes \iota)(\iota
\otimes T_1)$.  Therefore, the following should be no surprise.
\begin{prop} \label{dieT}
The map $T$ can be uniquely extended to a contraction from $A \oh
A \oh A$ to $A \ot A \ot A$. The extended $T$ is injective.
\end{prop}
\begin{proof}
First observe that $T(p \otimes q \otimes s) \in A \otimes A
\otimes A$ as e.g.\ $\phi(q)(1 \otimes r) \in A \otimes A$ and $(p \otimes 1)
\phi (s) \in A \otimes A$ (for all $p,q,r,s \in A$).
Let $(p_i)$, $(q_{ij})$ and $(r_j)$ be as before. Then by
lemma~\ref{41b} we have
$$\Bigl\| \sum_{i,j=1}^{n,m} (p_i \ot 1 \ot 1) \phi^{(2)}(q_{ij}) (1 \ot 1
\ot r_j) \Bigr\| \leq \|(p_i)\| \, \| (\phi^{(2)}(q_{ij})) \| \,
\|(r_j)\| \leq \|(p_i)\| \, \| (q_{ij}) \| \,
\|(r_j)\|.$$
So $T$ is a contraction for the Haagerup norm. Now suppose $x \in
A \oh A \oh A$ and $T(x)=0$. Let $\om \in A^*$ and apply $(\om \ot
\io \ot \io)$. This formally gives $T_1((\om \ot \io)T_2 \ot
\io)(x)=0$. To give a meaning to this we have to verify that the
map $(\om \ot \io) T_2$ from $A \oh A$ to $A$ is completely
bounded. But this is true because $T_2$ is completely bounded from
$A \oh A$ to $A \ot A$ and $(\om \ot \io)$ from $A \ot A$ to $A$.
Then the formula
$$ (\om \ot \io \ot \io)T(y) = T_1((\om \ot \io)T_2 \ot \io)(y)$$
is easily checked for $y \in A \od A \od A$ and then extended by
continuity. So we get indeed $((\om \ot \io)T_2 \ot \io)(x) = 0$
for all $\om \in A^*$ because $T_1$ is injective. But then
$$(\om \ot \io) T_2 ((\io \ot \io \ot \mu)(x)) = 0$$
for all $\om,\mu \in A^*$. Because $T_2$ is injective we get $(\io
\ot \io \ot \mu)(x)=0$ for all $\mu \in A^*$ and hence $x=0$.
\end{proof}

Now, we will use this triple Haagerup tensor product mainly to prove that the
left and right counit and antipode, as obtained in the previous section
coincide on a certain subalgebra, contained in $A_0 \cap A^\prime_0$.  It will
be shown in the examples that this subalgebra is still dense in most known
cases.  We will also obtain some more formulas for the counit and antipode.

We start with the following notation.
\begin{nota} \label{45}
Let $\ronda$ denote the space of elements $a \in A$
such that $b \otimes a \otimes c$ belongs to the range of $T$ on $A \otimes_h A
\otimes_h A$ for all $b,c \in A$.
\end{nota}

So when $a \in \ronda$ and $b,c \in A$ we can find $x \in A \oh A
\oh A$ such that $b \ot a \ot c = T(x)$. If now $\om \in A^*$ we
see, using the formula in the proof of proposition~\ref{dieT}, that
$$\om(b) a \ot c = (\om \ot \io \ot \io)T(x) = T_1((\om \ot
\io)T_2 \ot \io)(x).$$
So we get that $a \ot c$ belongs to the range of $T_1$ on $A \oh
A$. This is not sufficient for $a$ to be in $A_0$. In order to
have this one would need $T_1^{-1}(a \ot c) = y(1 \ot c)$ for some
element $y \in A \omh A$ and all $c \in A$, and this involves an
extra condition. If we assume however that $a \in A_0$ then we get
the following result.
\begin{prop} \label{46}
Let $a \in \ronda \cap A_0$. Let $b,c \in A$ and $x \in A \oh A
\oh A$ such that $b \ot a \ot c = T(x)$. Then
\begin{align*}
\eps(a)(b \ot c) &= E(x) &\text{with} \;  E(p \ot q \ot r) & = (p
\ot 1) \phi(q) (1 \ot r). \\
\intertext{For the antipode we get the following formulas.}
(b \ot 1 \ot 1)(1 \ot S(a) \ot 1)(1 \ot \phi(c)) &=F_1(x)
&\text{with} \;  F_1(p \ot q \ot r) &= (p \ot 1 \ot
1)\phi_{13}(q)\phi_{23}(r) \\
(\phi(b) \ot 1)(1 \ot S(a) \ot 1)(1 \ot \phi(c)) &= F_2(x)
&\text{with} \; F_2(p \ot q \ot r) &= (\phi(p) \ot 1)
\phi^{(2)}(q) (1 \ot \phi(r)) \\
bS(a) \ot c &= F_3(x) &\text{with} \; F_3(p \ot q \ot r) & = p \ot
qr \\
bS(a)c & = F_4(x) &\text{with} \;  F_4(p \ot q \ot r) & = pqr
\end{align*}
\end{prop}
Remark that before the proposition makes sense, we have to show
that $E,F_1,F_2,F_3$ and $F_4$ are all bounded for the Haagerup
norm, with the range in $A$ for $F_4$, in $A \ot A$ for $E$ and $F_3$ and in $A
\ot A \ot A$ for $F_1$ and $F_2$. This can be done completely
analogously as the extension of $T$ (see proposition~\ref{dieT}).

The formulas in the proposition can be easily illustrated for Hopf
algebras. So let $(H,\Delta,\eps,S)$ be a Hopf algebra and $x \in
H$. Then we clearly have, using Sweedler notation,
$$1 \ot x \ot 1 = T \Bigl( \sum S(x_{(1)}) \ot x_{(2)} \ot
S(x_{(3)}) \Bigr).$$
To illustrate the second formula for the antipode, we compute:
\begin{align*}
\sum \Delta_{12}(S(x_{(1)})) \Delta^{(2)}(x_{(2)})
\Delta_{23}(S(x_{(3)})) &= \sum S(x_{(2)}) x_{(3)} \ot S(x_{(1)})
x_{(4)} S(x_{(7)}) \ot x_{(5)} S(x_{(6)}) \\
&= \sum 1 \ot \eps(x_{(2)}) S(x_{(1)}) x_{(3)} S(x_{(5)})
\eps(x_{(4)}) \ot 1 \\
&= \sum 1 \ot S(x_{(1)}) x_{(2)} S(x_{(3)}) \ot 1 \\
&= 1 \ot \eps(x_{(1)}) S(x_{(2)}) \ot 1 = 1 \ot S(x) \ot 1.
\end{align*}
Analogously one can compute the other formulas.
\begin{proof}
Before this proposition we already remarked that
\begin{equation} \label{weereenvgl}
\om(b) a \ot c = T_1((\om \ot \io)T_2 \ot \io)(x).
\end{equation}
Because $a \in A_0$ we get from this that
$$\om(b) \eps(a) c = m((\om \ot \io)T_2 \ot \io) (x).$$
Now for $p,q,r \in A$ we have
$$m((\om \ot \io)T_2 \ot \io)(p \ot q \ot r) = (\om \ot \io)E(p
\ot q \ot r)$$
and so by continuity we have $\eps(a)b \ot c = E(x)$. This proves
the formula for the counit.

Now starting again with formula~\ref{weereenvgl} and using the
definition of the antipode (see~\ref{39}) we get
$$\om(b)(S(a) \ot 1) \phi(c) = \widetilde{T}((\om \ot \io)T_2 \ot
\io)(x).$$
But for all $p,q,r \in A$
$$\widetilde{T}((\om \ot \io)T_2 \ot \io)(p \ot q \ot r) = (\om
\ot \io \ot \io) F_1(p \ot q \ot r)$$
and so we get
$$(b \ot 1 \ot 1)(1 \ot S(a) \ot 1)(1 \ot \phi(c)) = F_1(x).$$
This gives the first formula for the antipode. Now choose $\om \in
A^*$ and apply $(\io \ot \io \ot \om)$ to this equation. This
gives
$$b \ot S(a) (\io \ot \om)\phi(c) = (\io \ot \io \ot \om)F_1(x).$$
But for $p,q,r \in A$ we have
$$(\io \ot \io \ot \om)F_1(p \ot q \ot r) = (p \ot 1) ((\io \ot
\io \ot \om)(\phi_{13}(q) \phi_{23}(r))).$$
In proposition~\ref{34} we have seen that
$$q \ot r \mapsto (\io \ot \io \ot \om)(\phi_{13}(q)
\phi_{23}(r))$$
extends to a bounded map from $A \oh A$ to itself. A slight
modification of that proof gives that this map is even completely
bounded. But this allows us to conclude that the map
$$p \ot q \ot r \mapsto (p \ot 1) ((\io \ot
\io \ot \om)(\phi_{13}(q) \phi_{23}(r)))$$
can be extended to a bounded map from $A \oh A \oh A$ to $A \oh
A$. We will denote this extension by $G$. So we get in $A \oh A$
the equality
\begin{equation} \label{dieG}
b \ot S(a)(\io \ot \om)\phi(c) = G(x).
\end{equation}
Then we may conclude that
\begin{equation} \label{weeralG}
\phi(b) (1 \ot S(a)(\io \ot \om)\phi(c)) = T_1(G(x)).
\end{equation}
Observe now that
\begin{align*}
T_1(G(p \ot q \ot r)) & = \phi(p) T_1 ((\io \ot \io \ot
\om)(\phi_{13}(q) \phi_{23}(r))) \\
&= \phi(p) (\io \ot \io \ot \om)(\phi^{(2)}(q) \phi_{23}(r)) \\
\intertext{as we already observed in the proof of
proposition~\ref{36}}
&= (\io \ot \io \ot \om)F_2(p \ot q \ot r)
\end{align*}
By continuity we conclude from formula~\ref{weeralG} that
$$\phi_{12}(b)(1 \ot S(a) \ot 1) \phi_{23}(c) = F_2(x).$$
This is the second formula for the antipode we had to prove. Next
we observe that
\begin{align*}
m(G(p \ot q \ot r)) & = p m((\io \ot \io \ot \om)(\phi_{13}(q)
\phi_{23}(r))) \\
& = p (\io \ot \om)\phi(qr) \\
\intertext{as we already observed in the proof of
proposition~\ref{38}}
&= (\io \ot \om)T_2(p \ot qr)
\end{align*}
Now we have to remark that $F_3$ is a contraction from $A \oh A
\oh A$ to $A \oh A$ and this allows us to write
$$m(G(x)) = (\io \ot \om)(T_2(F_3(x))).$$
Combined with formula~\ref{dieG} we get
$$T_2(bS(a) \ot c) = T_2(F_3(x))$$
and so $bS(a) \ot c = F_3(x)$. We already observed that $F_3(x)
\in A \oh A$ and so we get $bS(a)c = m(F_3(x))$. Clearly
$m(F_3(x))= F_4(x)$ and this concludes the proof of the
proposition.
\end{proof}

Now, some of these formulas are nicely symmetric.  This is the case e.g.\ for
the formula with $\varepsilon$ and for the second and the fourth formula involving
$S$.  Using the notation of remark~\ref{remark} this implies the following.

\begin{prop} \label{47}
If $a \in \ronda \cap A_0 \cap A^\prime_0$, then
$\varepsilon(a) = \varepsilon^\prime(a)$ and $S(a) = S^\prime(a)$.
\end{prop}

\begin{proof}
In a completely similar fashion we get for elements $a \in \ronda \cap A'_0$ formulas for $\varepsilon^\prime(a)$ and
$S^\prime(a)$ as in the previous proposition.  From these formulas, it
immediately follows that $\varepsilon(a) = \varepsilon^\prime(a)$ and $S(a) =
S^\prime(a)$.
\end{proof}

One might guess from the formulas in~\ref{46} that it is possible to define
$\varepsilon$ and $S$ on all of $\ronda$.  This is indeed the case.  It requires
some more arguments.  But it is not very relevant for the theory.  In some
sense, the space $\ronda$ is too big.  Moreover, there seems to be no
argument to show that $\ronda$ is an algebra which would be a natural
requirement.

Now, also the set $\ronda \cap A_0 \cap A^\prime_0$ will probably not be
an algebra.  Of course, as $A_0$ and $A_0^\prime$ are algebras, the
subalgebra of $A$
generated by $\ronda \cap A_0 \cap A^\prime_0$ will be a subalgebra of $A_0$
and of $A^\prime_0$ where $\varepsilon = \varepsilon^\prime$ and $S = S^\prime$.

From the theoretical point it would certainly be nice to have a more precise
understanding of all these sets.  On the other hand in practice it is just
important to have enough elements. This is what we will show in
the examples in the next section.

\section{Special cases and examples}
In this section, we will consider some special cases and examples.
Essentially, there are always three steps to be taken.  First there is the
construction of the comultiplication on the C$^*$-algebra.  Then one must
try to prove the injectivity of the maps $T_1$ and $T_2$ so that we do have a
Hopf C$^*$-algebra.  Finally one has to find the counit and antipode, as
constructed in section 3 (one-sided) and section 4 (two-sided) and hopefully
show that the respective domains are dense.
\subsection{Locally compact groups}
The main two examples to start with are $C_0(G)$ and $C^\ast_r(G)$ for a
locally compact group $G$.

The first result is obvious.
\begin{prop} \label{51}
Let $A$ be the C$^*$-algebra $C_0(G)$ of continuous
complex functions tending to $0$ at infinity on  a locally compact group $G$.
Let $\phi$ be defined as before by $\phi(f)(p,q) = f(pq)$ for $f \in C_0(G)$
and $p,q \in G$.  The $\phi$ is a comultiplication on $A$ and $(A, \phi)$ is a
Hopf C$^*$-algebra.
\end{prop}
The injectivity of the maps $T_1$ and $T_2$ follows as $A \otimes_h A$ is
contained in $A \otimes A$ and here $T_1$ and $T_2$ have continuous extensions
to $A \otimes A$ and these extensions are even C$^*$-automorphisms.

We give the following result mainly for reasons of presentation.  We will
formulate and prove more general results in proposition~\ref{56} and \ref{elinA}.

\begin{prop} \label{52}
Let $A = C_0(G)$ and $\phi$ as before.  Consider a
left Haar measure on $G$, take two functions $f$ and $g$ in $K(G)$, the space
of continuous complex functions with compact support in $G$ and consider $h \in
K(G)$ defined by
$$h(p) = \int f(pq)g(q)dq.$$
Then $h \in A_0$ (as defined in~\ref{31}) and $\eps(h) = h(e)$ and $S(h)(p) =
h(p^{-1})$.
\end{prop}

We will not give the proof here in detail as we will come back to this later.
But we can define a function $k$ of two variables by
$$k(p,q) = \int f(pr)g(qr)dr.$$
This is a bounded continuous function and one can show that it belongs to $A \otimes_{mh} A$.
If we apply $T_1$ to it, we obtain
\begin{align*}
(T_1 k)(p,q) & = k(pq,q) = \int f(pqr)g(qr)dr \\
& = \int f(pr)g(r)dr = h(p).
\end{align*}
It follows that $T_1(k) = h \ot 1$ and so $h \in A_0$.  Then $\eps(h) = k(p,p) = \int f(pr)g(pr)dr =
\int f(r)g(r)dr = h(e)$ for all $p$.  Moreover
\begin{align*}
k(q,pq) & = \int f(qr)g(pqr)dr \\
& = \int f(r)g(pr)dr \\
& = \int f(p^{-1}r)g(r)dr = h(p^{-1})
\end{align*}
and this illustrates the formula $S(h)(p) = h(p^{-1})$.

The next case is $C^\ast_r(G)$.  Here, the first step is well-known.

\begin{prop} \label{53}
Let $G$ be a locally compact group.  The reduced
C$^*$-algebra $C^\ast_r(G)$ of $G$ is the C$^*$-algebra generated by the
left convolution operators on $L^2(G)$, where $G$ is considered with the left
Haar measure.  If we define $W$ on $L^2(G) \otimes L^2(G)$ by
$$(W \xi)(p,q) = \xi(p,pq)$$
then $\phi(x) = W^*(x \otimes 1)W$ defines a comultiplication on
$C^\ast_r(G)$ and $(C^\ast_r(G),\phi)$ is a Hopf C$^*$-algebra.
\end{prop}
\begin{proof}
That $W$ induces a comultiplicaton on $C^\ast_r(G)$ is proved
in e.g.~\cite{Val}, proposition~4.2.
The only thing that remains to be shown is that the maps $T_1$ and
$T_2$ are injective on $C^\ast_r(G) \otimes_h C^\ast_r(G)$.  We will prove it
for $T_2$.  Then the injectivity of $T_1$ can be obtained easily as
$\phi(x) = \sigma \phi(x)$ (where $\sigma$ is the flip) so that
\begin{align*}
T_1(a \otimes b) & = \phi(a)(1 \otimes b) \\
& = ((1 \otimes b^\ast) \phi (a^\ast))^\ast \\
& = (\sigma(T_2(b^\ast \otimes a^\ast)))^\ast
\end{align*}
and we have seen that $a \otimes b \rightarrow b^\ast \otimes a^\ast$ is an
isometry of the Haagerup tensor product to itself.

Now remark that for $p,q \in C^*_r(G)$, $\om \in B(L^2(G))_*$ and
any right convolution operator $y$ on $L^2(G)$ we have
$$y \bigl((\io \ot \om)(T_2(p \ot q) W^*) \bigr) = y \bigl((\io \ot \om)((p \ot 1) W^*
(q \ot 1))\bigr) = p y \bigl((\io \ot \om)(W^*)\bigr) q$$
where we used the obvious fact that left and right convolution
operators commute. It is clear that for any $z \in B(L^2(G))$ we
can define a bounded map $M_z$ from $C^*_r(G) \oh C^*_r(G)$ to
$B(L^2(G))$, with norm majorized by $\|z\|$ and such that $M_z(p
\ot q) = p z q$. If now $x \in C^*_r(G) \oh C^*_r(G)$ and $T_2(x)=0$
we conclude by continuity from the computation above that
$M_z(x)=0$ for all $z$ of the form $y ((\io \ot \om)(W^*))$. But the
space $\{(\io \ot \om)(W^*) \mid \om \in B(L^2(G))_* \}$ is norm
dense in the space of multiplication operators by functions of
$C_0(G)$. So the linear span of the elements $y ((\io \ot
\om)(W^*))$ with $\om \in B(L^2(G))_*$ and $y$ a right convolution
operator, is norm dense in the compact operators on $L^2(G)$. So
we can conclude that $M_k(x)=0$ for all compact operators $k$ on
$L^2(G)$. In particular we get, for all $\xi_1,\xi_2,\eta_1,\eta_2
\in L^2(G)$
$$0 = \langle M_{\theta_{\xi_1} \theta_{\eta_2}^*} (x) \eta_1,
\xi_2 \rangle = (\om_{\xi_1,\xi_2} \ot \om_{\eta_1,\eta_2})(x)$$
where we used the notation $\theta_\xi$ for the rank one operator
from $\C$ to $L^2(G)$ given by $\theta_\xi(\lambda) = \lambda
\xi$.
From this it follows that $x=0$ and that completes the proof.
\end{proof}

Also here we can get enough elements in $A_0$.  Again this follows from a more
general result that we will prove later (see propositions~\ref{56} and \ref{elinA}).
But let us once more illustrate the theory by looking closer to this special case.
We have the following result.

\begin{prop} \label{54}
Let $f \in L^1(G)$ and let $\pi(f)$ denote left
convolution with $f$ on $L^2(G)$, then $\pi(f) \in A_0$ and $\eps(\pi(f)) =
\int f(p)dp$ and $S(\pi(f)) = \pi(\widetilde f)$ where $\widetilde{f}(p) =
\Delta(p)^{-1} f(p^{-1})$ and $\Delta$ is the modular function of
$G$.
\end{prop}

The argument here is essentially the following.  Denote by $\lambda_p$ the left
translation operator on $L^2$ by $p$ given by $(\lambda_p \xi)(q) =
\xi(p^{-1}q)$.  Put $h = \int f(p) (\lambda_p \otimes \lambda_{p^{-1}})dp$.
This will be in $C^\ast_r (G) \otimes_{mh} C^\ast_r(G)$ and

\begin{align*}
T_1(h) & = \int f(p)(\lambda_p \otimes  \lambda_p \lambda_{p^{-1}}) dp \\
& = \Bigl( \int f(p) \lambda_p dp \Bigr) \otimes 1.
\end{align*}
So $\int f(p) \lambda_p dp \in A_0$ and $\eps(\int f(p) \lambda_p dp) 1 =
\int f(p) \lambda_p \lambda_{p^{-1}} dp = \int f(p)dp  1$ and

\begin{equation*}
 \int f(p) (\lambda_{p^{-1}} \otimes \lambda_p \lambda_{p^{-1}})dp
= \int f(p) \lambda_{p^{-1}} dp = \int f(p^{-1}) \Delta (p^{-1}) \lambda_p dp.
\end{equation*}

For these two basic examples it is also possible to give elements in the space
$\ronda \cap A_0 \cap A^\prime_0$ (cf. previous section).  This is not so
relevant here, see propositions~\ref{58} and \ref{514}.

\subsection{Hopf C$^*$-algebras with corepresentations}
Now we will pass to a much more general situation and we will show how
corepresentations give rise to elements in the domain of the counit and the
antipode.

First recall the notion of a corepresentation as we will use it here.  In this
section, we will fix a Hopf C$^*$-algebra $(A, \phi)$.

\begin{defin} \label{55}
By a unitary corepresentation of $(A, \phi)$ on a
Hilbert space $H$ we mean a unitary element $u$ in $M(A \otimes
B_0(H))$, where $B_0(H)$ stands for the compact operators on
$H$, satisfying $(\phi \otimes \iota)(u) = u_{13} u_{23}$.
\end{defin}

We again use the leg numbering notation here : $u_{23}$ is $1 \otimes u$ in
$M(A \otimes A \otimes B_0(H))$ while $u_{13}$ is obtained from
$1 \otimes u$ by applying the flip on the first two factors.

\begin{prop} \label{56}
Let $u$ be a unitary corepresentation such that
$(\iota \otimes \omega)(u) \in A$ for all $\omega \in B(H)_\ast$.  Then
$(\iota \otimes \omega)(u) \in A_0$ for all $\omega \in B(H)_\ast$ and
\begin{equation*}
\eps (\iota \otimes \omega)(u) = \omega(1) \quad \quad
S(\iota \otimes \omega)(u) = (\iota \otimes \omega)(u^*)
\end{equation*}
\end{prop}
\begin{proof}
First remark that it is sufficient to look at the case $\om =
\om_{\xi,\eta}$ with $\xi,\eta \in H$. Because $\om \in B_0(H)^*$
we can take a representation $\pi$ of $B_0(H)$ on another Hilbert
space $K$ such that $\om = \om_{\xi,\eta} \circ \pi$ with $\xi,\eta \in
K$. Then we replace $u$ by $(\io \ot \pi)(u)$. So we may suppose
$\om = \om_{\xi,\eta}$ with $\xi,\eta \in H$. Let $(\rho_i)_{i \in I}$ be an
orthonormal basis for $H$ and define
$$p_i = (\io \ot \om_{\rho_i,\eta})(u) \quad\quad q_i = (\io \ot
\om_{\xi,\rho_i})(u^*).$$
We have $p_i,q_i \in A$ for all $i$ and by a slight modification
of lemma~\ref{strict} (or by lemma~9.5 in \cite{KV2}) we get that
$\sum_{i \in I} p_i p_i^*$ converges strictly to $(\io \ot
\om_{\eta,\eta})(u u^*)$. Also $\sum_{i \in I} q_i^* q_i$ is a
bounded strictly converging net. So by proposition~\ref{113} we
can define an element $x \in A \omh A$ by
$$x = \sum_{i \in I} p_i \ot q_i,$$
the sum being strictly convergent. Then we have
\begin{align*}
T_1(x) &=
\sum_{i \in I} \phi(p_{i})(1 \otimes q_{i}) \\ &= (\iota \otimes \iota \otimes
\omega_{\xi,\eta})((\phi \otimes \iota)(u) u_{23}^*) \\
& = (\iota \otimes \iota \otimes \omega) (u_{13} u_{23} u^\ast_{23}) \\
& = (\iota \otimes \iota \otimes \omega)(u_{13}) \\
& = (\iota \otimes \omega)(u) \otimes 1.
\end{align*}
So indeed $(\io \ot \om)(u) \in A_0$. Next we have $\eps((\io \ot
\om)(u))1 = \sum_{i \in I} p_i q_i = (\io \ot \om)(u u^*) = \om(1)
1$ and
\begin{align*}
S((\io \ot \om)(u)) \otimes 1 & = \sum_{i \in I} (1 \otimes p_{i})\phi(q_{i}) \\
& = (\iota \otimes \iota \otimes \omega_{\xi,\eta})(u_{23}(\phi \otimes
\iota)(u^\ast)) \\
& = (\iota \otimes \iota \otimes \omega)(u_{23} u^\ast_{23} u_{13}^*) \\
& = (\iota \otimes \iota \otimes \omega) (u^\ast_{13}) \\
& = (\iota \otimes \omega)(u^\ast) \otimes 1
\end{align*}
so that $S((\io \ot \om)(u)) = (\iota \otimes \omega)(u^\ast)$.
This completes the proof.
\end{proof}

We would like to comment a little on this result and the proof of it.

Observe that the left regular representation of a
locally compact quantum group will satisfy the requirement that $(\iota \otimes
\omega)(u) \in A$ for all $\omega$.  Remark that the trivial corepresentation $(1
\in M(A))$ does not have this property when $A$ has no identity.

Formally, we can write
\begin{equation*}
(\eps \otimes \iota)(u) = 1 \quad \quad
(S \otimes \iota)(u) = u^\ast
\end{equation*}

If $u$ and $v$ are corepresentations as in definition~\ref{55}, respectively on
$H$ and $K$, we can define a corepresentation $w$ on $H \otimes K$ by $w = u_{12}v_{13}$.
When $\omega_1 \in B(H)_\ast$
and $\omega_2 \in B(K)_\ast$ we find
$$(\iota \otimes \omega_1 \otimes \omega_2)(w) = (\iota \otimes
\omega_1)(u) (\iota \otimes \omega_2)(v)$$
and we see that we obtain a subalgebra of elements in $A_0$.  It is easy to
verify that $\eps$ is a homomorphism and $S$ a antihomomorphism on this
subalgebra.  Also $S(a)^\ast = (\iota \otimes \bar \omega)(u)$ when $a = (\iota
\otimes \omega)(u)$ so that again $S(a)^\ast$ is of this form.  And the formula
$S(S(a)^\ast)^\ast = a$ is obvious here.

Observe that the same formula for the antipode is obtained
in~\cite{Wor5}, Theorem~1.6, in the framework of manageable
multiplicative unitaries. Also for locally compact quantum groups
we have such a formula, but that is not so surprising because any
locally compact quantum group is a Hopf C$^*$-algebra. See
proposition~\ref{quantum} and the remark following it for a more
precise result. In the framework of Hopf algebras (or multiplier
Hopf algebras) these formulas are again well known, see \S\ 11.1.1
of~\cite{K-S} and \cite{VDZ}. So far about the left counit and antipode.

In a completely similar way, we can work with unitary elements $v$ in
$M(B_0(H) \otimes A)$ satisfying $(\iota \otimes \phi)(v) = v_{12} v_{13}$.
Then $(\om \otimes \iota)(v) \in A'_0$ and $\eps^\prime((\omega \otimes \iota)(v))
= \omega(1)$ and $S^\prime((\omega \otimes \iota)(v)) = (\omega \otimes
\iota)(v^\ast)$, when $(\om \ot \io)(v) \in A$ for all $\om$.
Now, because we can simply go from one type corepresentation to
the other by applying the flip, we get the following result.

\begin{prop} \label{57}
As before, let $u$ be a unitary corepresentation such
that $(\iota \otimes \omega)(u) \in A$ for all $\omega \in B(H)_*$.  Then $(\iota
\otimes \omega)(u) \in A_0 \cap A^\prime_0$ and $\eps$ and $\eps'$,
as well as $S$ and $S^\prime$ coincide on these elements.
\end{prop}
\begin{proof}
Let $v = \sigma u$ where $\sigma$ denotes the flip from $M(A
\otimes B_0(H))$ to $M(B_0 (H) \otimes A)$.  Then $(\iota \otimes
\phi)(v) = v_{12} v_{13}$.  If $a = (\omega \otimes \iota)(v)$ and $x = (\omega \otimes
\iota \otimes \iota)(v^\ast_{12} v_{13})$ one verifies, as in the proof of~\ref{56},
that $x \in A \otimes_{mh}A$ and that
\begin{align*}
T_2 (x) & = (\omega \otimes \iota \otimes \iota)(v^\ast_{12} (\iota \otimes
\phi)(v)) \\
& = (\omega \otimes \iota \otimes \iota)(v^\ast_{12} v_{12} v_{13}) \\
& = (\omega \otimes \iota \otimes \iota)(v_{13}) = 1 \otimes a.
\end{align*}
Hence $\eps^\prime(a)1 = (\omega \otimes \iota)(v^\ast v) = \omega (1)1$.  Also
\begin{align*}
1 \otimes S^\prime(a) & = (\omega \otimes \iota \otimes \iota)((\iota \otimes
\phi)(v^\ast)v_{12}) \\
& = (\omega \otimes \iota \otimes \iota)(v^\ast_{13} v^\ast_{12} v_{12})
\\
& = (\omega \otimes \iota \otimes \iota)(v^\ast_{13}) = 1 \otimes (\omega
\otimes \iota)(v^\ast).
\end{align*}
Together with the previous proposition this gives the result.
\end{proof}
As we have remarked before, taking into account all such corepresentations, we
get a subalgebra of $A_0 \cap A^\prime_0$ where $\eps = \eps^\prime$
and $S = S^\prime$.

Also remark that here we need $v^\ast v = 1$ and hence $u^\ast u = 1$ whereas
in the previous proposition we needed $uu^\ast  = 1$.

Now we come to the two-sided approach.  Here we can prove the following.

\begin{prop} \label{58}
Let $u$ and $w$ be unitary corepresentations on $H$ and $K$ respectively with
the property that $(\io \ot \om)(u)$ and $(\io \ot \om')(w)$ are in $A$ for $\omega \in B(H)_*$
and $\omega^\prime \in B(K)_\ast$ respectively.  Then the product of two
such elements belongs to $\ronda \cap A_0 \cap A^\prime_0$, as considered
in~\ref{45}.
\end{prop}

\begin{proof}
Take $\om \in B(H)_*$ and $\om' \in B(K)_*$. Just as in the proof
of proposition~\ref{56} we may suppose that $\om = \om_{\xi,\eta}$
and $\om'=\om_{\xi',\eta'}$. Write $a = (\io \ot \om)(u)$ and
$b=(\io \ot \om')(w)$. Use the flip map $\sigma$ and denote $v =
\sigma w \in M(B_0(K) \ot A)$. Then $b = (\om' \ot \io)(v)$ and
$$ab = (\om' \ot \io \ot \om)(u_{23} v_{12}).$$
From the previous proposition it follows that $a, b \in A_0 \cap
A'_0$, and because both are subalgebras of $A$ we get $ab \in A_0
\cap A'_0$. To show that $ab \in \ronda$ we consider
$$x = (\om_{\xi',\eta'} \ot \io \ot \io \ot \io \ot
\om_{\xi,\eta})(v_{12}^* u_{25}^* u_{35} v_{13} v_{14}^*
u_{45}^*).$$
We claim that for all $c,d \in A$ we have
$$(c \ot 1 \ot 1) x (1 \ot 1 \ot d) \in A \oh A \oh A$$
and
$$T((c \ot 1 \ot 1) x (1 \ot 1 \ot d)) = c \ot ab \ot d.$$
Let $(\rho_i)_{i \in I}$ be an orthonormal basis for $H$ and
$(\rho'_j)_{j \in J}$ for $K$. Define
\begin{align*}
p^\prime_k &= (\omega_{\rho^\prime_k,\eta^\prime} \otimes \iota)(v^\ast) &
p_i &= (\iota \otimes \omega_{\rho_i,\eta})(u^\ast) \\
q^\prime_{kl} &= (\omega_{\rho^\prime_l, \rho^\prime_k} \otimes \iota)(v) &
q_{ij} &= (\iota \otimes \omega_{\rho_j, \rho_i})(u) \\
r^\prime_l &= (\omega_{\xi^\prime,\rho^\prime_l} \otimes \iota)(v^\ast) & r_j &=
(\iota \otimes \omega_{\xi, \rho_j})(u^\ast).
\end{align*}
As in the proof of proposition~\ref{56} we get that the nets
$\sum_{k \in J} p'_k {p'}_k^*$ and $\sum_{i \in I} p_i p_i^*$ are both bounded
and strictly converging. So also
$$\sum_{(i,k) \in I \times J} p'_k p_i p_i^* {p'}_k^*$$
will be bounded and strictly converging. Analogously
$$\sum_{(j,l) \in I \times J} r_j^* {r'}_l^* r_l r_j$$
is bounded and strictly converging. It is clear that $(q_{ij}
q'_{kl})_{(i,k),(j,l)}$ is a bounded (infinite) $I \times J$ by $I
\times J$ matrix over $A$. Then it is not hard to see that for
every $c,d \in A$ the net
$$\sum_{i,k,j,l} c p'_k p_i \ot q_{ij} q'_{kl} \ot r'_l r_j d$$
is norm convergent in $A \oh A \oh A$ (see also the estimate in
the proof of proposition~\ref{514}). The limit in $A \ot A \ot A$
is given by $(c \ot 1 \ot 1)x(1 \ot 1 \ot d)$. Applying $T$ to
this element gives, after some computations,
\begin{equation} \label{lang}
(c \ot 1 \ot 1) (\omega_{\xi',\eta'} \otimes \iota \otimes \iota \otimes \iota \otimes
\omega_{\xi,\eta})(v^\ast_{12} u^\ast_{25} u_{25} u_{35} u_{45} v_{12} v_{13} v_{14}
v^\ast_{14} u^\ast_{45})(1 \ot 1 \ot d)
\end{equation}
where we used
$$(\iota \otimes \phi^{(2)} \otimes \iota)(1 \otimes u) = u_{25} u_{35} u_{45}
\quad \quad
(\iota \otimes \phi^{(2)} \otimes \iota)(v \otimes 1) = v_{12} v_{13}
v_{14}.$$
But now observe that the element in formula~\ref{lang} equals
\begin{align*}
&(c \ot 1 \ot 1) (\omega^\prime \otimes \iota \otimes \iota \otimes \iota \otimes \omega)
(v^\ast_{12} u_{35} u_{45} v_{12} v_{13} u^\ast_{45}) (1 \ot 1 \ot d)\\
& \hspace{3cm} = (c \ot 1 \ot 1)(\omega^\prime \otimes \iota \otimes \iota \otimes \iota \otimes
\omega)(u_{35}u_{45} u^\ast_{45} v^\ast_{12} v_{12} v_{13}) (1 \ot 1 \ot d)\\
& \hspace{3cm} = (c \ot 1 \ot 1) (\omega^\prime  \otimes \iota \otimes \iota \otimes \iota \otimes
\omega)(u_{35} v_{13})(1 \ot 1 \ot d) \\
& \hspace{3cm} = c \otimes (\omega^\prime \otimes \iota \otimes \omega)(u_{23} v_{12})
\otimes d \\
& \hspace{3cm} = c \otimes ab \otimes d.
\end{align*}
So we get that $c \ot ab \ot d$ belongs to the range of $T$ on $A
\oh A \oh A$ for all $c,d \in A$. This gives $ab \in \ronda$.
\end{proof}

Remark that here we need both $u^\ast u = 1$ and $uu^\ast = 1$, and similarly
for $w$.

Let us illustrate some of the formulas that we obtained in section
4.

If we apply the map $p \otimes q \otimes r \mapsto (p \otimes 1)\phi(q)(1
\otimes r)$ to the element $x$ in the proof of the proposition, we get

\begin{align*}
(\omega^\prime \otimes \iota \otimes \iota \otimes \omega)(v^\ast_{12}
u^\ast_{24} u_{24} u_{34} v_{12} v_{13} v^\ast_{13} v^\ast_{34})
& = (\omega^\prime \otimes \iota \otimes \iota \otimes \omega)(v^\ast_{12}
u_{34} v_{12} u^\ast_{34}) \\
& = (\omega^\prime \otimes \iota \otimes \iota \otimes \omega)(v^\ast_{12}
v_{12} u_{34} u^\ast_{34})  \\ & = \omega^\prime (1)  \omega (1) 1 \otimes 1
\end{align*}
and this illustrates the formula
$$\eps (a) (b \otimes c) = \sum (p_i \otimes 1)\phi (q_{ij})(1 \otimes
r_j)$$
we obtained in proposition~\ref{46}.

If we apply $p \otimes q \otimes r \rightarrow (\phi(p) \otimes 1)
\phi^{(2)} (q) (1 \otimes \phi(r))$ to $x$ we get
\begin{align*}
& (\omega^\prime \otimes \iota \otimes \iota \otimes \iota
\otimes \omega) (v^\ast_{13} v^\ast_{12} u^\ast_{35} u^\ast_{25} u_{25} u_{35}
u_{45} v_{12} v_{13} v_{14} v^\ast_{14} v^\ast_{13} u^\ast_{45} u^\ast_{35}) \\
& \hspace{6cm} = (\omega^\prime \otimes \iota \otimes \iota \otimes \iota \otimes
\omega)(v^\ast_{13} v^\ast_{12} u_{45} v_{12} u^\ast_{45} u^\ast_{35}) \\
& \hspace{6cm} = (\omega^\prime \otimes \iota \ot \io \otimes \iota \otimes \omega)(v^\ast_{13}
u^\ast_{35}) \\
& \hspace{6cm} = 1 \otimes (\omega^\prime \otimes \iota)(v^\ast) (\iota \otimes
\omega)(u^\ast) \otimes 1 \\
& \hspace{6cm} = 1 \otimes S(b)S(a) \otimes 1 = 1 \otimes S(ab) \otimes 1
\end{align*}
and this illustrates the formula
$$(\phi(b) \otimes 1)(1 \otimes S(a) \otimes 1)(1 \otimes \phi(c)) = \sum(\phi
(p_i) \otimes 1) \phi^{(2)} (q_{ij}) (1 \otimes \phi (r_j)).$$

Finally, look at the map $p \otimes q \otimes r \rightarrow p \otimes qr$
and apply it to $x$.  This gives
\begin{align*}
(\omega^\prime \otimes \iota \otimes \iota \otimes \omega)(v^\ast_{12}
u^\ast_{24} u_{34} v_{13} v^\ast_{13} u^\ast_{34})
& = (\omega^\prime \otimes \iota \otimes \iota \otimes \omega)(v^\ast_{12}
u^\ast_{24}) \\
& = (\omega^\prime \otimes \iota)(v^\ast)(\iota \otimes \omega)(u^\ast) \otimes 1
\\
& = S(b)S(a) \otimes 1 = S(ab) \otimes 1.
\end{align*}
which illustrates the formula
$$bS(a) \otimes c = \sum p_i \otimes q_{ij}r_j.$$

So we can conclude that for Hopf C$^*$-algebras $(A, \phi)$ with enough
corepresentations $u$ satisfying $(\iota \otimes \omega) (u) \in A$ for all
$\omega \in B(H)_\ast$, we get a dense subalgebra $A_1$, namely the set
of such elements $(\iota \otimes u)(u)$, contained in $A_0 \cap A^\prime_0$.  By
taking $A^2_1$, we  even get a dense subalgebra contained in $\ronda \cap A_0
\cap A'_0$.  But already $\eps = \eps^\prime$ and $S = S^\prime$ on
$A_1$.

\subsection{Locally compact quantum groups}

The last class of Hopf C$^*$-algebras we will consider here are the ones
admitting faithful Haar measures.

Let us first concentrate on the problem of the injectivity of the maps $T_1$
and $T_2$.  Clearly we have the following result.

\begin{prop} \label{59}
Let $(A, \phi)$ be a pair of a C$^*$-algebra $A$
with identity and a comultiplication $\phi$.  Assume there is a faithful state
$h$ on $A$ which is both left and right invariant (i.e.\ $(\iota \otimes
h)\phi(a) = h(a)1$ and $(h \otimes \iota)\phi(a) = h(a)1$ for all $a \in A$),
then the maps $T_1$ and $T_2$ are injective on $A \otimes_h A$.
\end{prop}
\begin{proof}
Let $x \in A \otimes_h A$ and assume that $T_1(x)=0$. Now observe
that for all $a \in A$, $\om \in A^*$ and $p,q \in A$
$$(h \ot \om)(\phi(a) T_1(p \ot q)) = (h \ot \om) (\phi(ap) (1 \ot
q)) = h(a p \om(q)) = h(a (\io \ot \om)(p \ot q)).$$
By continuity we may conclude that
$$h(a (\io \ot \om)(x)) = (h \ot \om)(\phi(a) T_1(x)) = 0.$$
By the faithfulness of $h$ we get $(\io \ot \om)(x)=0$ for all
$\om \in A^*$. So $x=0$.
\end{proof}

This is the case of a compact quantum group with a faithful Haar measure.
Remark that the Hopf $^\ast$-algebra obtained as matrix elements of
finite-dimensional unitary corepresentations will be (properly) contained in
the algebra that we obtained previously in this section and that of course the
counit and antipode, when restricted to this Hopf $^\ast$-algebra gives the
counit and the antipode of the Hopf $^*$-algebra.

This proposition indicates the relation between the faithfulness of the Haar
measures and the injectivity of the maps $T_1$ and $T_2$.  Unfortunately, in
the general case, where we have to work with weights (instead of states), the
situation is much more complicated.

We will use chapter~1 of~\cite{KV2} as a reference for weights on
C$^*$-algebras and use freely the notations and results introduced
there. If $\vphi$ is a weight on a C$^*$-algebra $A$ we will use the
notation $\Nfi=\{\, x \in A \mid \vphi(x^*x) < \infty \, \}$ and
$\Mfi^+=\{ \, x \in A^+ \mid \vphi(x) < \infty \, \}$.
We use the abbreviation l.s.c. for lower semicontinuous.
As an analogue for the left invariance of a measure we
introduce the following:
\begin{defin} \label{510}
Let $A$ be a C$^*$-algebra and $\phi$ a
comultiplication on $A$.  A weight $\varphi$ on $A$ is called left invariant if
$\varphi((\omega \otimes \iota)\phi(a)) = \varphi(a)\omega(1)$ whenever $a \in
{\cal M}^+_\varphi$ and $\omega \in A^\ast_+$.
\end{defin}

Once we have a left invariant weight on a Hopf C$^*$-algebra
$(A,\phi)$ we can construct elements in $A_0$. The result we
obtain is closely related to proposition~\ref{56}: with this left
invariant weight one can also construct the left regular
corepresentation $W$ and then the elements in $A_0$ given in the
following proposition are of the form $(\io \ot \om)(W)$.
Nevertheless we cannot simply apply proposition~\ref{56} because
with the given assumptions it seems impossible to prove that $W$
is a unitary (and not only a partial isometry) and that $W$
belongs to $M(A \ot B_0(H))$.
\begin{prop} \label{elinA}
Let $(A,\phi)$ be a Hopf C$^*$-algebra and suppose that $\vphi$ is
a l.s.c. densely defined left invariant weight on
$A$. Then for every $a,b \in \Nfi$ we have
$$(\io \ot \vphi)(\phi(a^*)(1 \ot b)) \in A_0$$
and
$$S((\io \ot \vphi)(\phi(a^*)(1 \ot b))) = (\io \ot \vphi)((1 \ot
a^*) \phi(b)) \quad\quad \eps((\io \ot \vphi)(\phi(a^*)(1 \ot b))) =
\vphi(a^*b).$$
When moreover $\{(\io \ot \om) \phi(a) \mid \om \in
A^*, a \in A\}$ spans a dense subspace of $A$ and when $\vphi$ is
faithful, then $A_0$ is dense in $A$.
\end{prop}
\begin{proof}
Let $(H_\vphi,\pifi,\lafi)$ be a GNS-construction for $\vphi$. Let
$(\xi_i)_{i \in I}$ be an orthonormal basis for $H_\vphi$ and
define for $\eta \in H_\vphi$ the operator
$\theta_\eta(\lambda)=\lambda \eta$ for $\lambda \in \C$. Now we
will use the map $\io \ot \lafi$ introduced in proposition~1.26 of
\cite{KV2}. This map takes values in $\cL(A, A \ot \Hfi)$, the
space of adjointable maps between the Hilbert $A$-modules
$A$ and $A \ot \Hfi$. Of course, for $a \in M(A)$ and $b \in \Nfi$
we have $(\io \ot \lafi)(a \ot b) (x) = ax \ot \lafi(b)$ for all $x
\in A$, but it is possible to extend $\io \ot \lafi$ to a much
larger domain, see~\cite{KV2}.
Then one has $(\io \ot \lafi)(y)^* (\io \ot \lafi)(x) = (\io \ot
\vphi)(y^*x)$ for all $x$ and $y$ in the domain of $\io \ot
\lafi$.
With the help of this we define
$$p_i^* = (1 \ot \theta_{\xi_i}^*)(\io \ot \lafi) \phi(a)
\quad\text{and}\quad q_i = (1 \ot \theta_{\xi_i}^*)(\io \ot
\lafi)\phi(b).$$
Then $p_i,q_i \in \cL(A)= M(A)$ for all $i \in I$. But if we fix
$i \in I$ we can approximate $q_i$ in norm by elements of the form
$$(1 \ot \theta_{\lafi(x)}^*)(\io \ot \lafi) \phi(b)$$
with $x \in \Nfi$. This equals $(\io \ot \vphi)((1 \ot x^*)
\phi(b))$ and hence belongs to $A$ by 2.5 in \cite{KV2}. So
$p_i,q_i \in A$ for all $i \in I$. Now consider the net
$$\sum_{i \in I} p_i p_i^* = ((\io \ot \lafi)\phi(a))^* \Bigl(1 \ot
\sum_{i \in I} \theta_{\xi_i} \theta_{\xi_i}^* \Bigr) (\io \ot
\lafi)\phi(a).$$
This net is strictly convergent to
$$((\io \ot \lafi)\phi(a))^* (\io \ot \lafi) \phi(a) = (\io \ot
\vphi)\phi(a^*a) = \vphi(a^* a) 1.$$
Analogously $\sum_{i \in I} q_i^* q_i$ is bounded and strictly
convergent. So by proposition~\ref{113} we can define, with
strict convergence
$$x = \sum_{i \in I} p_i \ot q_i \in A \omh A.$$
Then we get
\begin{align*}
T_1(x) &= \sum_{i \in I} \phi(p_i) (1 \ot q_i) \\
&= \sum_{i \in I} (\io \ot \io \ot \lafi)(\phi^{(2)}(a))^* (1 \ot
1 \ot \theta_{\xi_i} \theta_{\xi_i}^*) (\io \ot \io \ot \lafi)(1
\ot \phi(b)) \\
&= (\io \ot \io \ot \vphi)(\phi^{(2)}(a^*)(1 \ot \phi(b))) \\
&= (\io \ot \vphi)(\phi(a^*)(1 \ot b)) \ot 1
\end{align*}
by left invariance. All these manipulations are very rigorous when
interpreted in the sense of~\cite{KV2}. So we may conclude that $(\io \ot \vphi)(\phi(a^*)(1 \ot
b)) \in A_0$. Then we can make analogous computations as above to get
the stated formulas for $S$ and $\eps$.

Assume now the extra conditions stated in the proposition. Let
$\mu \in A^*$ such that $\mu(a)=0$ for all $a \in A_0$. By the
Hahn-Banach theorem it is enough to prove $\mu = 0$. By the result
above we have $\vphi((\mu \ot \io)\phi(a^*) b) = 0$ for all $a,b
\in \Nfi$. Because $\vphi$ is faithful we conclude that $(\mu \ot
\io)\phi(a^*) = 0$ for all $a \in \Nfi$ and hence for all $a \in
A$ by density and continuity. Then we get $\mu((\io \ot
\om)\phi(a)) = 0$ for all $a \in A$ and $\om \in A^*$.
By the assumed density we get $\mu = 0$.
\end{proof}

The density of $A_0$ implies that the range of $T_1$ on $A \otimes_h A$ is
dense in $A \otimes A$.  And this implies that $\phi(A)(1 \otimes A)$ is dense
in $A \otimes A$.  On the other hand, the density of $\phi(A)(1 \otimes A)$ in
$A \otimes A$ will imply that $\{(\iota \otimes \omega)\phi(a) | \omega \in
A^\ast, a \in A \}$ will be dense in $A \otimes A$.  So we find that under the
given circumstances, these two density conditions are in some sense equivalent.

A similar result holds when there is a right invariant weight and this will
give the density of $A^\prime_0$.

If moreover the left or the right regular representation is unitary, it will
follow from~\ref{58} that also $\ronda \cap A_0 \cap A^\prime_0$ is dense in  $A$.
This is the case when $(A, \phi)$ is a locally compact quantum group in the
sense of~\cite{KV1} and \cite{KV2}, see proposition~\ref{quantum}.  Roughly speaking, a locally compact quantum group is a pair
$(A,\phi)$ of a C$^*$-algebra $A$ with a comultiplication $\phi$ such that
the spaces $\phi(A)(1 \otimes A)$ and $\phi(A)(A \otimes 1)$ are dense in $A
\otimes A$ and such that there exist left and right invariant faithful l.s.c.
weights that are KMS.

However, to obtain the density of $\ronda \cap A_0 \cap A^\prime_0$ in $A$,
one can do with weaker assumptions. We have the following result.
\begin{prop} \label{514}
Let $(A,\phi)$ be a Hopf C$^*$-algebra.  Assume
there is a densely defined left invariant l.s.c. weight $\varphi$ and a densely
defined right invariant l.s.c. weight $\psi$.  Then, for all $a,b \in {\cal
N}_\varphi$ and $c,d \in {\cal N}_\psi$ we have that
$$x := (\psi \otimes \iota \otimes \varphi)((c^\ast \otimes 1 \otimes
1)\phi^{(2)} (a^\ast d)(1 \otimes 1 \otimes b))$$
belongs to $\ronda \cap A_0 \cap A^\prime_0$ and $\eps(x) = \psi (c^\ast
d)\varphi(a^\ast b)$ and
\begin{align*}
S(x) & =  (\psi \otimes \iota \otimes \varphi)((c^\ast \otimes 1 \otimes
1)\phi_{13} (a^\ast d) \phi_{23}(b)) \\
& = (\psi \otimes \iota \otimes \varphi) (\phi_{12} (c^\ast) \phi_{13} (a^\ast
d)(1 \otimes 1 \otimes b)) \\
& = (\psi \otimes \iota \otimes \varphi) (\phi_{12}(c^*) \phi^{(2)} (a^\ast d)
\phi_{23} (b)).
\end{align*}
\end{prop}
\begin{proof}
Again, as in the proof of the previous proposition, we will use
the results of chapter~1 of~\cite{KV2}. Implicitly a lot of weight
theory techniques that appeared in that paper will be used. Also
remark that a more elementary but very long proof for
this result is given in~\cite{V}.
So we can consider
$$(\lapsi \ot \io \ot \io)(\phi^{(2)}(d)) \in \cL(A \ot A, \Hpsi
\ot A \ot A).$$
If now $T \in \cL(A \ot A, \Hpsi \ot A \ot A)$ it is possible to
define
$$(\io \ot \io \ot \pifi)(T) \in \cL(A \ot \Hfi, \Hpsi \ot A \ot
\Hfi)$$
such that
$$(\io \ot \io \ot \pifi)(T) (a \ot \pifi(b)\xi)=(\io \ot \io \ot
\pifi( \cdot) \xi)(T(a \ot b))$$
where
$$(\io \ot \io \ot \pifi( \cdot) \xi)(\eta \ot x \ot y) = \eta \ot
x \ot \pifi(y) \xi.$$
This should not be too surprising. If $\eta \in \Hpsi$ we have on
the one hand
$$(\io \ot \pifi)((\theta_\eta^* \ot 1 \ot 1) T) \in M(A \ot
B_0(\Hfi))$$
and on the other hand
$$(\theta_\eta^* \ot 1 \ot 1) (\io \ot \io \ot \pifi)(T) \in \cL(A
\ot \Hfi)$$
and this agrees with the usual identification of $M(A \ot
B_0(\Hfi))$ and $\cL(A \ot \Hfi)$.
So we can consider the element $P \in \cL(A \ot \Hfi, \Hpsi \ot
A)$ given by
$$P= [ (\pipsi \ot \io \ot \io)(\io \ot \io \ot
\lafi)(\phi^{(2)}(a)) ]^* [(\io \ot \io \ot \pifi)(\lapsi \ot \io
\ot \io)(\phi^{(2)}(d)) ].$$
Now let $(e_j)_{j \in J}$ and $(f_i)_{i \in I}$ be orthonormal
bases for $\Hfi$ and $\Hpsi$ respectively. Define for $i \in I$
and $j \in J$:
$$q_{ij}=(\theta_{f_i}^* \ot 1) P (1 \ot \theta_{e_j}).$$
Then we have that $q_{ij} \in \cL(A)= M(A)$. But for $x \in \Nfi$
and $y \in \Npsi$ we have
\begin{align*}
&(\theta_{\lapsi(y)}^* \ot 1) P (1 \ot \lafi(x)) \\
& \quad = \bigl[ (\pipsi \ot \io \ot \io)(\io \ot \io \ot
\lafi)(\phi^{(2)}(a)) (\theta_{\lapsi(y)} \ot 1) \bigr]^* \\
& \hspace{5cm} \bigl[ (\io \ot \io \ot \pifi) (\lapsi \ot \io \ot \io)
(\phi^{(2)}(d)) (1 \ot \theta_{\lafi(x)}) \bigr] \\
& \quad = (\lapsi \ot \io \ot \lafi)(\phi^{(2)}(a)(y \ot 1 \ot
1))^* (\lapsi \ot \io \ot \lafi)( \phi^{(2)}(d)(1 \ot 1 \ot x)) \\
& \quad = (\psi \ot \io \ot \vphi)((y^* \ot 1 \ot 1) \phi^{(2)}(a^*d) (1 \ot 1 \ot x)) \\
& \quad = (\psi \ot \io) \bigl( (y^* \ot 1) \phi((\io \ot
\vphi)(\phi(a^*d)(1 \ot x))) \bigr).
\end{align*}
By the results of \cite{KV2} one has $(\io \ot
\vphi)(\phi(a^*d)(1 \ot x)) \in \Npsi$ and so the expression above
belongs to $A$. Then we get $q_{ij} \in A$ for all $i$ and $j$.
From its definition it follows now immediately that $(q_{ij})_{i
\in I,j \in J}$ is a bounded matrix over $A$. Define also
\begin{equation*}
p_i^* = (\theta_{f_i}^* \ot 1) (\lapsi \ot \io)(\phi(c))
\quad\quad
r_j = (1 \ot \theta_{e_j}^*)(\io \ot \lafi)(\phi(b)).
\end{equation*}
As in the proof of the previous proposition we get $\sum_{i \in I}
p_i p_i^*$ and $\sum_{j \in J} r_j^* r_j$ strictly converging and
bounded. But then we have for any $u, v \in A$ that
$$\sum_{i \in I, j \in J} u p_i \ot q_{ij} \ot r_j v$$
converges in norm in $A \oh A \oh A$, because for $I_0 \subseteq
I$ and $J_0 \subseteq J$ finite we have
$$\Bigl\| \sum_{i \in I_0, j \in J_0} u p_i \ot q_{ij} \ot r_j v
\Bigr\|_h \leq \Bigl\| \sum_{i \in I_0} u p_i p_i^* u^*
\Bigr\|^{1/2} \, \|(q_{ij})_{i \in I,j \in J} \| \, \Bigl\| \sum_{j \in J_0} v^* r_j^*
r_j v \Bigr\|^{1/2}.$$
In order to apply $T$ to this we have to realize that (use the
abbreviation $\iota^4$ for $\io \ot \io \ot \io \ot \io$)
$$\phi^{(2)}(q_{ij}) = (\theta_{f_i}^* \ot 1 \ot 1 \ot 1) \bigl[
(\pipsi \ot \io^4)(\io^4 \ot \lafi)(\phi^{(4)}(a)) \bigr]^*
\bigl[ (\io^4 \ot \pifi)(\lapsi \ot \io^4)(\phi^{(4)}(d)) \bigr]
(1 \ot 1 \ot 1 \ot \theta_{e_j}^*).$$
Then it is not too hard to see that with strict convergence
\begin{align*}
&\sum_{i \in I,j \in J} (p_i \ot 1 \ot 1) \phi^{(2)}(q_{ij})(1 \ot
1 \ot r_j) \\
&\quad = \bigl[ (\pipsi \ot \io^4) \bigl((\io^4 \ot
\lafi)(\phi^{(4)}(a)) \bigr) (\lapsi \ot \io^3)(\phi_{12}(c))
\bigr]^* \\
& \hspace{5cm}
\bigl[ (\io^4 \ot \pifi) \bigl( (\lapsi \ot \io^4)(\phi^{(4)}(d)) \bigr) (\io^3
\ot \lafi)(\phi_{34}(b)) \bigr] \\
&\quad = (\lapsi \ot \io \ot \io \ot \io \ot \lafi)(\phi^{(4)}(a)
\phi_{12}(c))^* (\lapsi \ot \io \ot \io \ot \io \ot
\lafi)(\phi^{(4)}(d) \phi_{45}(b)) \\
&\quad = (\psi \ot \io \ot \io \ot \io \ot \vphi)(\phi_{12}(c^*)
\phi^{(4)}(a^*d) \phi_{45}(b)) \\
&\quad = (\psi \ot \io \ot \io \ot \io \ot \vphi) (\phi \ot \io \ot
\phi) ((c^* \ot 1 \ot 1) \phi^{(2)}(a^*d) (1 \ot 1 \ot b)) \\
&\quad = 1 \ot (\psi \ot \io \ot \vphi)((c^* \ot 1 \ot 1)
\phi^{(2)}(a^*d) (1 \ot 1 \ot b)) \ot 1 \\
&\quad = 1 \ot x \ot 1.
\end{align*}
So we get indeed that $x \in \ronda$ and for all $u,v \in A$ we
have
$$u \ot x \ot v = T \Bigl( \sum_{i,j} u p_i \ot q_{ij} \ot r_j z
\Bigr).$$
Because we can write $x = (\io \ot \vphi)(\phi((\psi \ot \io)((c^*
\ot 1) \phi(a^*d))) (1 \ot b))$ and because $(\psi \ot \io)((c^*
\ot 1) \phi(a^*d)) \in \Nfi^*$ we get $x \in A_0$ by the previous
proposition. Analogously $x \in A'_0$.
The formulas to be proven for $S(x)$ can be obtained from
proposition~\ref{46} with a similar calculation as above.
\end{proof}

We finish this section by the following result:
\begin{prop} \label{quantum}
If $(A, \phi)$ is a (reduced) locally compact quantum group in the sense of~\cite{KV1} and \cite{KV2},
then $(A,\phi)$ is a Hopf C$^*$-algebra and the space $\ronda \cap A_0 \cap
A^\prime_0$ is a core for the antipode obtained in the theory of locally
compact quantum groups.
\end{prop}
So, in particular the discrete and compact quantum groups will be
Hopf C$^*$-algebras and also the C$^*$-Kac algebras will be.
\begin{proof}
We use of course the notations and conventions from \cite{KV1} and
\cite{KV2}. So let $(A,\phi)$ be a reduced locally compact quantum group. We let
$A$ act on the Hilbert space $H$, which is the GNS-space of the left Haar
weight $\vphi$. We first proof that $(A,\phi)$ is a Hopf C$^*$-algebra. We
will work along the same lines as in the proof of proposition~\ref{53}.
Let $W \in B(H \ot H)$ be the multiplicative unitary associated with
$(A,\phi)$ and observe that for any $\om \in B(H)_*$ and any $y \in A'$ we
have
$$y (\om \ot \io)(W T_1(p \ot q)) = y (\om \ot \io) ((1 \ot p)W (1 \ot q))
= p y (\om \ot \io)(W) q.$$
We can define for any $z \in B(H)$ a bounded map $M_z$ from $A \oh A$
to $B(H)$ such that
$$M_z(p \ot q) = p z q$$
for all $p,q \in A$. Then it is clear that for any $x \in A \oh A$ the map
$$z \mapsto M_z(x)$$
is strongly continuous on bounded subsets of $B(H)$. Let now $x \in A \oh
A$ and suppose $T_1(x)=0$. Then it follows from the computation above that
$M_z(x)=0$ for all $z$ of the form $y (\om \ot \io)(W)$, with $y \in A'$
and $\om \in B(H)_*$. Now we claim two things.
First, the closed linear span of $\{ \, y a \mid y \in JAJ \, , \, a \in \hat{A} \, \}$
is a C$^*$-algebra that acts non-degenerately on $H$.
Secondly this C$^*$-algebra is strongly dense in $B(H)$.
Here $J$ denotes the modular conjugation of the weight $\vphi$, and so
$JAJ \subseteq A'$. After these two claims have been proven we may
conclude that $M_z(x)=0$ for all $z \in B(H)$ and exactly as in the proof
of proposition~\ref{53} this is enough to conclude $x=0$.

To prove the first claim denote the stated closed linear space by $L$. It
is then enough to prove that for all $\om,\mu \in B(H)_*$
\begin{equation} \label{laatste}
(\mu \ot \io)(W) J (\io \ot \om)(W) J \in L.
\end{equation}
But this expression equals
$$J(\rho \ot \io)(W^*) (\io \ot \om)(W) J$$
where $\rho(x) = \overline{\mu(\hat{J} x \hat{J})}$ and $\hat{J}$ is the
modular conjugation of the dual weight $\hat{\vphi}$. Here we used the
equality $W (\hat{J} \ot J) = (\hat{J} \ot J) W^*$. Now observe that
\begin{align*}
(\rho \ot \io)(W^*) (\io \ot \om)(W) &= (\rho \ot \io \ot \om)(W^*_{12}
W_{23}) \\
&= (\rho \ot \io \ot \om)(W_{13} W_{23} W_{12}^*) \\
\intertext{where we used the pentagon equation $W_{12} W_{13} W_{23} =
W_{23} W_{12}$. So we get}
(\rho \ot \io)(W^*) (\io \ot \om)(W) &\in \text{Closed span} \{ (\rho \ot \io \ot \om) (W_{23} W_{12}^*) \mid
\rho,\om \in B(H)_* \} \\
&= \text{Closed span} \{ x y \mid x \in A, y \in \hat{A} \}.
\end{align*}
Hence we can conclude that the expression in~\ref{laatste} belongs to the
closed linear span of
$$\{JxyJ \mid x \in A, y \in \hat{A} \} = \{JxJ \, JyJ \mid x \in A,y \in
\hat{A} \} = L$$
because $J \hat{A} J = \hat{A}$. So we have proved the first claim.

To prove the second claim first observe that $A^{\prime \prime} \cap {\hat{A}}' = \C
1$. If $x \in A^{\prime \prime} \cap {\hat{A}}'$ then
$$\widetilde{\phi}(x) = W^*(1 \ot x) W = 1 \ot x$$
because $W \in M(A \ot \hat{A})$, and with $\widetilde{\phi}$ the extension of
$\phi$ to the von Neumann algebra $A^{\prime \prime}$. But this implies $x \in \C 1$.
By taking the commutant we get $(J A J \cup \hat{A})^{\prime \prime} = B(H)$ and from
this follows the second claim. So we have proven the injectivity of
$T_1$.

The injectivity of $T_2$ follows easily by using the unitary
antipode $R$. Just observe that for $p,q \in A$
$$\chi (R \ot R) (T_2(p \ot q)) = T_1( \chi (R \ot R) (p \ot q))$$
where $\chi$ denotes the flip map extended to $M(A \ot A)$. Because $R$ is
an anti-automorphism of $A$ we can extend $\chi (R \ot R)$ to an
isomorphism of $A \oh A$ (compare to proposition~\ref{112}). Then the
injectivity of $T_2$ follows from the already proven injectivity of $T_1$.

From proposition~8.3 in \cite{KV2} it follows that $\{ (\io \ot \om)(W) \mid \om \in B(H)_* \}$ is a
core for the antipode of the locally compact quantum group $(A,\phi)$, which we
shall denote by $S_1$.
Because $S_1$ is an anti-homomorphism also the linear span of the
products of two such elements will give a core. It is not hard to
finetune the proofs of propositions~5.33 and 5.45 in~\cite{KV2} to obtain
that $A_0$ is a subspace of the domain of $S_1$. All these remarks,
together with proposition~\ref{58} give that $\ronda \cap A_0 \cap A'_0$
is a core for $S_1$.
\end{proof}

Finally remark that we have, in case of a locally compact quantum
group, that both $A_0$ and $A'_0$ are subspaces of the domain of
$S_1$ and $S_1$ coincides with $S$ and $S'$ on respectively $A_0$
and $A'_0$. In particular $S$ and $S'$ coincide on $A_0 \cap
A'_0$, although in the general situation we could only prove this
on $\ronda \cap A_0 \cap A'_0$.

\subsection{Further examples}
Combining the results of~\cite{Wor2} and \cite{Baa4} one gets that
Woronowicz' quantum $E(2)$ group is a locally compact quantum
group in the sense of~\cite{KV1} and \cite{KV2} and hence is a
Hopf C$^*$-algebra. Also the quantum Heisenberg group will be a
Hopf C$^*$-algebra because it is even a C$^*$-Kac algebra.

Finally the Kac systems of~\cite{B-S} give rise to Hopf
C$^*$-algebras and this can be proved completely analogously to
proposition~\ref{quantum}. This provides us with nice examples as
the quantum double of a locally compact group and the quantum
Lorentz group as studied by Podle\'s en Woronowicz in~\cite{PW}.


\begin{thebibliography}{VD}
\bibitem{Abe} {\sc E. Abe},
Hopf algebras. {\it Cambridge Tracts in Mathematics}, {\bf 74} {\it
Cambridge University Press, Cambridge} (1980).

\bibitem{B-S}  {\sc S. Baaj  \&  G. Skandalis},
Unitaires multiplicatifs et dualit\'e pour les produits crois\'es de
C$^*$-alg\`ebres. {\it Ann. scient. \'{E}c. Norm. Sup., $4{}^e$
s\'{e}rie}, {\bf 26} (1993), 425--488.

\bibitem{Baa4} {\sc S. Baaj}, Repr\'{e}sentation r\'{e}guli\`{e}re du groupe quantique des d\'{e}placements de Woronowicz.
{\it Ast\'erisque} {\bf 232} (1995), 11--48.

\bibitem{B-P} {\sc D.P. Blecher \& V.I. Paulsen},
Tensor Products of Operator Spaces. {\it J. Funct. Anal.} {\bf 99}
(1991), 262--292.

\bibitem{B-P2} {\sc D.P. Blecher}, Tensor Products of Operator
Spaces II. {\it Can. J. Math.} {\bf 44 (1)} (1992), 75--90.

\bibitem{EKR} {\sc E.G. Effros, J. Kraus \& Z.-J. Ruan}, On two
quantized tensor products. In {\it Operator algebras, mathematical physics,
and low dimensional topology.}, ed. R. Herman et al., A.K.
Peters. Res. Notes Math., Boston, Mass. {\bf 5} (1993), 125--145.

\bibitem{E-R} {\sc E.G. Effros \& Z.-J. Ruan}, Self-Duality for
the Haagerup Tensor Product and Hilbert Space Factorizations. {\it
J. Funct. Anal.} {\bf 100} (1991), 257--284.

\bibitem{E}  {\sc  M. Enock \&  J.-M. Schwartz},
Kac Algebras and Duality of Locally Compact Groups.
{\it Springer-Verlag, Berlin}  (1992).

\bibitem{E-V} {\sc  M. Enock \& J.-M. Vallin},
C$^*$-alg\`ebres de Kac et alg\`ebres de Kac.
{\it Proc. London Math. Soc. (3)} {\bf 66} (1993), 619--650.

\bibitem{K-S} {\sc A. Klimyk \& K. Schm\"udgen}, Quantum Groups
and Their Representations. {\it Springer-Verlag, Berlin} (1997).

\bibitem{KV1} {\sc J. Kustermans \& S. Vaes}, A simple
definition  for locally compact quantum groups. {\it C.R. Acad. Sci., Paris, S\'er. I}, {\bf 328 (10)}
(1999), 871--876.

\bibitem{KV2} {\sc J. Kustermans \& S. Vaes}, Locally compact
quantum groups. {\it Preprint University College
Cork, KU Leuven} (1999).

\bibitem{Mas-Nak} {\sc T. Masuda \& Y. Nakagami}, A von Neumann algebra framework for the duality of the quantum groups. {\it Publ. RIMS, Kyoto University} {\bf 30} (1994), 799--850.

\bibitem{Ng} {\sc C.-K. Ng}, Coactions and crossed products of
Hopf C$^*$-algebras. {\it Proc. London Math. Soc.} {\bf (3) 72}
(1996), 638--656.

\bibitem{Ped}   {\sc  G.K. Pedersen}, C$^*$-algebras and their automorphism groups. {\it Academic Press,  London} (1979).

\bibitem{JPP} {\sc J.-P. Pier}, Amenable locally compact groups.
{\it Pure \& Applied mathematics, John Wiley, New York} (1984).

\bibitem{PW} {\sc P. Podle\'{s} \& S.L. Woronowicz}, Quantum deformation of Lorentz group. {\it Commun. Math. Phys.} {\bf 130} (1990), 381--431.


\bibitem{Tak}  {\sc  M. Takesaki},
Theory of Operator Algebras I.
{\it Springer-Verlag, New York} (1979).

\bibitem{V} {\sc S. Vaes}, Hopf-C$^*$-algebra's. {\it Masters thesis, KU Leuven} (1998).

\bibitem{Val} {\sc J.-M. Vallin}, C$^*$-alg\`{e}bres de Hopf et C$^*$-alg\`{e}bres de Kac. {\it Proc. London Math. Soc.} (3){\bf 50} (1985), 131--174.

\bibitem{VD2} {\sc A. Van Daele}, Multiplier Hopf Algebras. {\it Trans. Amer. Math. Soc.} {\bf 342} (1994), 917--932.

\bibitem{VD4} {\sc A. Van Daele}, An algebraic framework for group duality.
{\it Adv. in Math.} {\bf 140} (1998), 323--366.

\bibitem{VDZ} {\sc A. Van Daele \& Y. Zhang}, Corepresentation
theory of multiplier Hopf algebras I. Preprint KU Leuven (1997), to appear in {\it Int. J.
Math.}

\bibitem{Wor5}  {\sc  S.L. Woronowicz}, From multiplicative unitaries to quantum groups.
{\it Int. J. Math.} {\bf Vol.~7, No.~1} (1996), 127--149.

\bibitem{Wor2}  {\sc  S.L. Woronowicz}, Quantum $E(2)$ group and its Pontryagin dual.
{\it Lett. Math. Phys.} {\bf 23} (1991), 251--263.

\end{thebibliography}
\end{document}